\documentclass[12pt,reqno]{amsart}
\usepackage{amsmath,amsfonts,amsthm,amsopn,amssymb}
\usepackage{cite,marginnote}
\usepackage{color,enumitem,graphicx}
\usepackage[colorlinks=true,urlcolor=blue,
citecolor=red,linkcolor=blue,linktocpage,pdfpagelabels,
bookmarksnumbered,bookmarksopen]{hyperref}
\usepackage[english]{babel}

\usepackage[left=2.9cm,right=2.9cm,top=2.8cm,bottom=2.8cm]{geometry}

\numberwithin{equation}{section}

\makeindex

\newtheorem{Remark}{Remark}[section]

\newtheorem{lemma}{Lemma}[section]
\newtheorem{iteration lemma}{iteration Lemma}[section]

\newcommand{\s}{\section}

\newcommand{\R}{\mathbb R}

\newcommand{\bt}{\begin{theorem}}
\newcommand{\et}{\end{theorem}}
\newcommand{\bl}{\begin{lemma}}
\newcommand{\el}{\end{lemma}}
\newcommand{\bd}{\begin{definition}}
\newcommand{\ed}{\end{definition}}
\newcommand{\bc}{\begin{corollary}}
\newcommand{\ec}{\end{corollary}}
\newcommand{\bp}{\begin{proof}}
\newcommand{\ep}{\end{proof}}
\newcommand{\bx}{\begin{example}}
\newcommand{\ex}{\end{example}}
\newcommand{\bi}{\begin{exercise}}
\newcommand{\ei}{\end{exercise}}
\newcommand{\bo}{\begin{proposition}}
\newcommand{\eo}{\end{proposition}}
\newcommand{\br}{\begin{remark}}
\newcommand{\er}{\end{remark}}
\newcommand{\beq}{\begin{equation}}
\newcommand{\eeq}{\end{equation}}
\newcommand{\ba}{\begin{align}}
\newcommand{\ea}{\end{align}}
\newcommand{\bn}{\begin{enumerate}}
\newcommand{\en}{\end{enumerate}}
\newcommand{\bg}{\begin{align*}}
\newcommand{\bcs}{\begin{cases}}
\newcommand{\ecs}{\end{cases}}

\newcommand{\bean}{\begin{eqnarray*}}
\newcommand{\eean}{\end{eqnarray*}}

\def\R{\mathbb{R}}

\def\bd{\mathrm{bd}\,}

\title[Localized semiclassical states]{Localized semiclassical states for Hamiltonian elliptic systems in dimension two}

\author[H. Zhang]{Hui Zhang}
\author[M. B.\ Yang]{Minbo Yang}
\author[J. J. Zhang]{Jianjun Zhang}
\author[X.~X.~Zhong]{Xuexiu Zhong}

\address[H.\ Zhang]{\newline\indent Department of  Mathematics,
Jinling Institute of Technology,
\newline\indent
Nanjing 211169, PR China
\newline\indent and
\newline\indent Department of Mathematics,
Nanjing University,
\newline\indent
Nanjing 210093, PR China}
\email{\href{mailto:huihz0517@126.com}{huihz0517@126.com}}

\address[M. B.\ Yang]{\newline\indent Department of Mathematics
\newline\indent
Zhejiang Normal University
\newline\indent
Jinhua 321004, PR China}
\email{\href{mailto:mbyang@zjnu.edu.cn}{mbyang@zjnu.edu.cn}}

\address[J. J. \ Zhang]{\newline\indent College of Mathematics and Statistics, Chongqing Jiaotong University,
\newline\indent
Chongqing 400074, PR China}
\email{\href{mailto:zhangjianjun09@tsinghua.org.cn}{zhangjianjun09@tsinghua.org.cn}}

\address[X.~X.~Zhong]{\newline\indent South China Research Center for Applied Mathematics and Interdisciplinary Studies
\newline\indent
South China Normal University
\newline\indent
Guangzhou 510631, PR China}
\email{\href{mailto:zhongxuexiu1989@163.com}{zhongxuexiu1989@163.com}}

\thanks{(1) Corresponding author: Minbo Yang ({\tt mbyang@zjnu.edu.cn})}

\thanks{(2) Hui Zhang was supported by China Postdoctoral Science Foundation (No.2021M691527). Minbo Yang was supported by NSFC(No.11971436, No.12011530199) and ZJNSF(No.LZ22A010001, No.LD19A010001). Xuexiu Zhong was supported by the NSFC (No.11801581), Guangdong Basic and Applied Basic Research Foundation (2021A1515010034),Guangzhou Basic and Applied Basic Research Foundation(No.202102020225). Jianjun Zhang was supported by NSFC (No.11871123).}

\subjclass[2000]{35J20; 35B25; 35J61.}
\keywords{Hamiltonian elliptic system; Trudinger-Moser inequality; Penalization method; Concentration; Asymptotically linear.}

\begin{document}

\begin{abstract}
 In this paper, we consider the Hamiltonian elliptic system in dimension two\begin{equation}\label{1.5}\aligned
\left\{ \begin{array}{lll}
-\epsilon^2\Delta u+V(x)u=g(v)\ & \text{in}\quad \mathbb{R}^2,\\
-\epsilon^2\Delta v+V(x)v=f(u)\ & \text{in}\quad \mathbb{R}^2,
\end{array}\right.\endaligned
\end{equation}
where $V\in C(\mathbb{R}^2)$ has local minimum points, and $f,g\in C^1(\mathbb{R})$ are assumed to be either superlinear or asymptotically linear at infinity and of subcritical exponential growth in the sense of Trudinger-Moser inequality. Under only a local condition on $V$, we obtain a family of semiclassical states concentrating around local minimum points of $V$. In addition, in the case that $f$ and $g$ are superlinear at infinity,  the decay and positivity of semiclassical states are also given.
The proof is based on a reduction method, variational methods and penalization techniques.
\end{abstract}

\maketitle

\s{Introduction and main results}
\renewcommand{\theequation}{1.\arabic{equation}}
Consider the time-dependent system of coupled Schr\"{o}dinger equations
\begin{equation}\label{1.4}\aligned
\left\{ \begin{array}{lll}
i\hbar\frac{\partial \varphi}{\partial t}=-\frac{\hbar^2}{2m}\Delta\varphi+W(x)\varphi-\frac{\partial H(\varphi,\psi)}{\partial \psi}, \ & x\in\mathbb{R}^N,\ t>0,\\
i\hbar\frac{\partial \psi}{\partial t}=-\frac{\hbar^2}{2m}\Delta\psi+W(x)\psi-\frac{\partial H(\varphi,\psi)}{\partial \varphi}, \ & x\in\mathbb{R}^N,\ t>0,
\end{array}\right.\endaligned
\end{equation}
where $N\geq2$, $(\varphi,\psi)$ represents the wave function of the state of an electron, $i$ is the
imaginary unit, $m$ is the mass of a particle, $\hbar$ is the Planck constant, $W(x)$ is a continuous potential, and $H$ is a coupled nonlinear function modeling various types of interaction
effect among many particles. System (\ref{1.4}) arises in many fields of physics, especially in nonlinear optics and
Bose-Einstein condensates theory. One of the most interesting problem for system (\ref{1.4}) is to look for standing waves, i.e.
$$(\varphi(x,t),\psi(x,t))=\bigl(u(x)e^{-iE/\hbar t}, v(x)e^{-iE/\hbar t}\bigr).$$
Assume $$\frac{\partial H(\varphi,\psi)}{\partial \psi}=e^{-iE/\hbar t}\frac{\partial H(u,v)}{\partial v},\quad \frac{\partial H(\varphi,\psi)}{\partial \varphi}=e^{-iE/\hbar t}\frac{\partial H(u,v)}{\partial u},$$
then $(\varphi,\psi)$ satisfies (\ref{1.4}) if and only if $(u,v)$ is a solution of the system
\begin{equation}\label{1.3}\aligned
\left\{ \begin{array}{lll}
-\epsilon^2\Delta u+V(x)u=\frac{\partial H(u,v)}{\partial v}, \ & x\in\mathbb{R}^N,\\
-\epsilon^2\Delta v+V(x)v=\frac{\partial H(u,v)}{\partial u}, \ & x\in \mathbb{R}^N,
\end{array}\right.\endaligned
\end{equation}
where $\epsilon^2=\frac{\hbar^2}{2m}$ and $V(x)=W(x)-E$. The existence, multiplicity and concentration of solutions of system (\ref{1.3}) are widely investigated and the scenario changes remarkably from the higher
dimensional case $N\geq3$ to the planar case $N=2$. In particular, $N=2$ affects the notion of critical growth which is the maximal admissible growth for the nonlinearities to preserve the variational structure of the problem, for details see \cite{Cassani1,ZhangCPDE} for coupled Schr\"{o}dinger systems and see \cite{TaoLiuYang,AlvesYang,ZhangDO} for a single Schr\"{o}dinger equation for instance. In the case that $\epsilon=1$ in (\ref{1.3}), the results on the existence and properties of solutions can be seen in \cite{BD,DY,LY,BDR,BS,DING} for $N\geq3$ and \cite{DG1,Cassani1,Cassani2,BDRT} for $N=2$ for example.

In the present paper, we are interested in the existence and asymptotic behavior of solutions as $\epsilon\rightarrow0$, which is called semiclassical problem and  the associated solutions are called semiclassical states, which possesses an important physical interest in
describing the translation from quantum to classical mechanics.
The semiclassical states of system (\ref{1.3}) has been intensively studied and most of the results are focused on the case $N\geq3$.  When $N\geq3$,  \'{A}vila and Yang \cite{Avila} considered a system with zero Neumann boundary condition
 \begin{equation*}\aligned
\left\{ \begin{array}{lll}
-\epsilon^2\Delta u+u=|v|^{p-1}v\ \text{in}\quad \Omega,\\
-\epsilon^2\Delta v+v=|u|^{q-1}u\ \text{in}\quad \Omega,\\
\frac{\partial u}{\partial n}=\frac{\partial v}{\partial n}=0 \ \text{on}\ \partial\Omega,
\end{array}\right.\endaligned
\end{equation*}
where  $\Omega$ is a bounded domain in $\mathbb{R}^N$ with smooth boundary $\partial\Omega$, and $2<p,q<2N/(N-2)$, they obtained
a family of positive solutions which concentrate near $\partial\Omega$. For related results, see \cite{RamosYang, Pisto}. Alves et al. \cite{alves} dealt with the system in the whole space
\begin{equation}\label{1.10}\aligned
\left\{ \begin{array}{lll}
-\epsilon^2\Delta u+u=Q(x)|v|^{p-1}v\ & \text{in}\quad \mathbb{R}^N,\\
-\epsilon^2\Delta v+v=K(x)|u|^{q-1}u\ & \text{in}\quad \mathbb{R}^N,
\end{array}\right.\endaligned
\end{equation}
where $2<p,q<2N/(N-2)$, $Q$ and $K$ are positive bounded functions, they proved that system (\ref{1.10}) admits a family of solutions concentrating at a point $x_0\in\mathbb{R}^N$ where
a related functional realizes its minimum energy.
Subsequently, Ramos and  Soares \cite{Ramos1} treated the system
\begin{equation}\label{1.5}\aligned
\left\{ \begin{array}{lll}
-\epsilon^2\Delta u+V(x)u=g(v)\ & \text{in}\quad \mathbb{R}^N,\\
-\epsilon^2\Delta v+V(x)v=f(u)\ & \text{in}\quad \mathbb{R}^N,
\end{array}\right.\endaligned
\end{equation}
where $V\in C(\mathbb{R}^N)$ satisfies
\begin{equation} \label{1.6} 0<\min_{x\in\mathbb{R}^N}V(x)<\liminf_{|x|\rightarrow\infty}V(x)\leq+\infty,\end{equation}
and $f,g\in C^1(\mathbb{R})$ are power-type nonlinearities having superlinear and subcritical growth at infinity. They showed the existence of ground states of system (\ref{1.5}) that concentrating at global minimum points of $V$ as $\epsilon\rightarrow0$. 
Later,  Ramos and Tavares \cite{Ramos2} considered (\ref{1.5}) in an open domain $\Omega$ of $\mathbb{R}^N$, and replaced the global condition (\ref{1.6}) by assuming $V(x)\geq a>0$ is locally H\"{o}lder continuous, and
\begin{equation*}\aligned \inf_{\Lambda_i}V<\inf_{\partial\Lambda_i}V,\ \text{where}\ \Lambda_i\subset\Omega  \text{ are bounded and mutually disjoint}, i=1,\cdots,k.\endaligned\end{equation*}
By means of a reduction method
and a penalization technique, the authors in \cite{Ramos2} showed that system (\ref{1.5}) possesses positive solutions $u_\epsilon, v_\epsilon$ and both of them have $k$ local maximum points $x_{i,\epsilon}\in \Lambda_i$, $i=1,\cdots,k$, as $\epsilon\rightarrow0$. In \cite{Ding2}, Ding et al. studied the system with competing potentials
\begin{equation}\label{1.8}\aligned
\left\{ \begin{array}{lll}
-\epsilon^2\Delta u+u+V(x)v=W(x)f(|z|)v\ & \text{in}\quad \mathbb{R}^N,\\
-\epsilon^2\Delta v+v+V(x)u=W(x)f(|z|)u\ & \text{in}\quad \mathbb{R}^N,
\end{array}\right.\endaligned
\end{equation}
where $z=(u,v)$, $V,W\in C^1(\mathbb{R}^N)$ satisfies global conditions as (\ref{1.6}) and $f\in C^1(\mathbb{R}^+)$ satisfies Ambrosetti-Rabinowitz condition (AR condition for briefly). Applying the reduction method
and Nehari technique, the
authors in \cite{Ding2} obtained a family of
ground states of system (\ref{1.8}) that concentrating around a concrete set characterized by the minimum points of $V$ or the maximum points of $W$.

On the contrary, the related results of semiclassical states for system (\ref{1.3}) in the plane are few.
In $\mathbb{R}^2$, the natural growth restriction on the nonlinearity is
defined by

\noindent{\bf Definition 1.1.} A function $h:\mathbb{R}\rightarrow\mathbb{R}$ is of subcritical exponential growth  if $$\lim_{|t|\rightarrow+\infty}\frac{|h(t)|}{e^{\alpha t^2}}=0, \quad \forall\alpha>0,$$ and critical exponential growth if there exists $\alpha_0>0$
such that
 $$
\lim_{|t|\rightarrow+\infty}\frac{|h(t)|}{e^{\alpha t^2}}
=\left\{\begin{array}{lll}
0\ & \text{if}\  \alpha>\alpha_0,\\
+\infty\ & \text{if}\ \alpha<\alpha_0.
\end{array}\right.
$$
Moreover, Definition 1.1 is derived from Trudinger-Moser inequality as follows.
\begin{lemma}\label{l1.2} (\cite{Cao}) If $\alpha>0$ and $u\in H^1(\mathbb{R}^2)$, then
$\int_{\mathbb{R}^2}\bigl(e^{\alpha u^2}-1\bigl)dx<+\infty.$
Moreover, if $|\nabla u|^2_2\leq1$, $|u|_2\leq M<+\infty$, and $\alpha<\alpha_0=4\pi$, then there exists a constant $C$, which depends only on $M$ and $\alpha$, such that
\begin{equation*}\int_{\mathbb{R}^2}\bigl(e^{\alpha u^2}-1\bigl)dx<C(M,\alpha).\end{equation*}
\end{lemma}

In \cite{ZhangCPDE}, Cassani and Zhang dealt with the system in dimension two
\begin{equation}\label{1.1}\aligned
\left\{ \begin{array}{lll}
-\epsilon^2\Delta u+V(x)u=g(v)\ & \text{in}\quad \mathbb{R}^2,\\
-\epsilon^2\Delta v+V(x)v=f(u)\ & \text{in}\quad \mathbb{R}^2,
\end{array}\right.\endaligned
\end{equation}
where $f,g\in C(\mathbb{R})$ are superlinear at infinity and of critical exponential growth. Under the global condition (\ref{1.6}) they showed the
existence of ground states concentrating around
global minimum points of the potential $V$ by the method of generalized Nehari manifold. In addition, Cassani and Zhang \cite{ZhangCPDE} studied the positivity and decay of ground states by virtue of a priori estimates of solutions for the limit system. The question on whether there is a family of solutions of system (\ref{1.1}) concentrating around local minimum points remains open.
In addition, to the best of our knowledge, there are no relevant results concerned with the system (\ref{1.1}) with asymptotically linear nonlinearity even in higher dimension $N\geq3$. So another interesting problem is whether similar results can be obtained for system (\ref{1.1}) with asymptotically linear nonlinearity. The purpose of the present paper is to fill these gaps and give affirmative answers. The assumptions on $V$, $f$ and $g$ are given as follows.
\vskip 0.1 true cm
\noindent(V$_1$) $V\in C(\mathbb{R}^2)\cap L^\infty(\mathbb{R}^2)$ and  $\inf_{\mathbb{R}^2}V>0$.
\vskip 0.1 true cm
\noindent(V$_2$) There exists an open bounded set $\Lambda\subset\mathbb{R}^2$ with smooth boundary  such that
$$V_0:=\min_{\Lambda}V<\min_{\partial\Lambda}V.$$
\vskip 0.05 true cm
\noindent(H$_1$) $f,g\in C^1(\mathbb{R},\mathbb{R})$, $f(t)=g(t)=0$ if $t\leq0$, $f(t)=o(t)$ and $g(t)=o(t)$ as $t\rightarrow0^+$.\\
\noindent(H$_2$) $tf'(t)>f(t)\ \text{and}\  tg'(t)>g(t),$ $\forall t>0$.
\vskip 0.1 true cm
\noindent(H$_3$) (i) $g'$ and $f'$ are of subcritical exponential growth in the sense of Definition 1.1;\\
\indent\ \ (ii) there is $\theta>2$ such that $\theta F(t)\leq f(t)t$ and $\theta G(t)\leq g(t)t$ for any $t>0$.

We are also concerned with the asymptotically linear case and replace (H$_3$) by the following condition.
\vskip 0.1 true cm
\noindent(H$'_3$) (i)\ $|V|_{L^\infty(\mathbb{R}^2)}<\lim_{t\rightarrow+\infty}
\frac{f(t)}{t}=\lim_{t\rightarrow+\infty}
\frac{g(t)}{t}=l_0<+\infty$;\\
\indent\ \ (ii) $f(t)t-2F(t)\rightarrow+\infty$ and $g(t)t-2G(t)\rightarrow+\infty$ as $t\rightarrow+\infty$.

\begin{Remark}\label{r1.1}
\noindent(1) If (H$_1$) and (H$_2$) are satisfied, then
\begin{equation}\label{1.2}t^2f'(t)>f(t)t>2F(t)>0\ \text{and}\  t^2g'(t)>g(t)t>2G(t)>0,\ \ \quad\forall t>0.\end{equation}
\noindent(2) The assumption (H$'_3$)-(i) can be weaken to be  $|V|_{L^\infty(\mathbb{R}^2)}<\min\{l_1,l_2\}$, where $l_1=\lim_{t\rightarrow+\infty}
\frac{f(t)}{t}<+\infty$ and $l_2=\lim_{t\rightarrow+\infty}
\frac{g(t)}{t}<+\infty$.
\end{Remark}

Setting $$\mathcal{V}=\{x\in \Lambda: V(x)=V_0\},$$
we give the first  result.
\vskip 0.05 true cm
\textbf{Theorem 1.1.} {\it Suppose that (V$_1$), (V$_2$), (H$_1$), (H$_2$) and either (H$_3$) or (H$'_3$) are satisfied. Then for sufficiently small $\epsilon>0$,\\
\noindent(1) the equation (\ref{1.1})
has a nontrivial solution $z_\epsilon=(\varphi_\epsilon,\psi_\epsilon)$ in $H^1(\mathbb{R}^2)\times H^1(\mathbb{R}^2)$;\\
\noindent(2) If additionally $V$ is uniformly continuous, then \\
\noindent(i) there exist $x^1_\epsilon, x^2_\epsilon, x_\epsilon\in\mathbb{R}^2$ be  (global) maximum points of $|\varphi_\epsilon|$, $|\psi_\epsilon|$ and $|\varphi_\epsilon|+|\psi_\epsilon|$ respectively, such that
 $$\lim_{\epsilon\rightarrow0}dist(x_\epsilon,\mathcal{V})=0, \ \text{and}\ \lim_{\epsilon\rightarrow0}dist(x^i_\epsilon,\mathcal{V})=0, \ i=1,2.$$}
{\it \noindent(ii)  $\bigl(\varphi_\epsilon(\epsilon x+x_\epsilon),\psi_\epsilon(\epsilon x+x_\epsilon)\bigr)$ and $\bigl(\varphi_\epsilon(\epsilon x+x^i_\epsilon),\psi_\epsilon(\epsilon x+x^i_\epsilon)\bigr),\ i=1,2$, converge  in $H^1(\mathbb{R}^2)\times H^1(\mathbb{R}^2)$ to ground states of
\begin{equation}\label{1.1.1}\aligned
\left\{ \begin{array}{lll}
-\Delta u+V_0u=g(v)\ & \text{in}\quad \mathbb{R}^2,\\
-\Delta v+V_0v=f(u)\ & \text{in}\quad \mathbb{R}^2.
\end{array}\right.\endaligned
\end{equation}}

\begin{Remark}
As we said before, there are no relevant results concerned with the singularly perturbed Hamiltonian system (\ref{1.1}) with asymptotically linear nonlinearity even in higher dimension $N\geq3$. So Theorem 1.1 about asymptotically linear nonlinearity can be extended to the higher dimensional case $N\geq3$.
\end{Remark}

In order to investigate the sign and decay of
 solutions of system (\ref{1.1}), we require in addition the following condition.
\vskip 0.05 true cm
\noindent(H$_4$) There exist $p,q>1$ such that $f(t)\geq t^q$ and $g(t)\geq t^p$ for small $t>0$.
\vskip 0.05 true cm
\textbf{Theorem 1.2.} {\it Assume (V$_1$), (V$_2$) and (H$_1$)-(H$_3$) are satisfied and $V$ is uniformly continuous. If in addition (H$_4$) holds and replacing $f,g$ with their odd extensions, then $\varphi_\epsilon,\ \psi_\epsilon$, $x^1_\epsilon$, $x^2_\epsilon$ obtained in Theorem 1.1 satisfy\\
\noindent(1) $\varphi_\epsilon\psi_\epsilon>0$ in $\mathbb{R}^2$;\\
\noindent(2) $x^1_\epsilon$, $x^2_\epsilon$ are unique, and
$$\lim_{\epsilon\rightarrow0}|x^1_\epsilon-x^2_\epsilon|/\epsilon=0.$$
Moreover, for some $C, c>0$ one has
$$|\varphi_\epsilon(x)|\leq C exp\bigl(-\frac c\epsilon|x-x^1_{\epsilon}|\bigr),\ |\psi_\epsilon(x)|\leq C exp\bigl(-\frac c\epsilon|x-x^2_{\epsilon}|\bigr).$$}
\begin{Remark}\label{r1.2} We would like to point out that,  system (\ref{1.1}) is quiet different from (\ref{1.8}) even in the same dimension, since (\ref{1.8}) can be reduced to the following form
\begin{equation*} Lz+V(\epsilon x)z=W(\epsilon x)f(|z|)z, \ z=(u,v),\end{equation*}
where $$\aligned L=\left( \begin{array}{lll}\ \ 0 &-\Delta+1\\
-\Delta+1 &\ \ 0\end{array}\right).\endaligned$$
Similar forms see also the Reaction-diffusion system and Dirac equation in \cite{DINGXU1,DINGXU2}. Hence, we cannot argue trivially as \cite{DINGXU1,DINGXU2}.
\end{Remark}
\begin{Remark} Observe that the behavior of $V$ outside $\Lambda$ is irrelevant, so if we suppose the following condition\\
\noindent(V$'_2$) there exist mutually disjoint open bounded domains $\Lambda_j$ with smooth boundary $\partial\Lambda_j$, $j=1,2,\cdots,k$ and constants
$a_1<\cdots<a_k$ such that
$$a_j:=\min_{\Lambda_j}V<\min_{\partial\Lambda_j}V,$$
then for sufficiently small $\epsilon>0$,  system (\ref{1.1})
admits $k$ nontrivial solutions which have similar properties as in Theorems 1.1 and 1.2.
\end{Remark}

Our results are based on variational methods, reduction and penalization arguments. Compared with related existing issues  of system (\ref{1.1}) in \cite{ZhangCPDE}, the distinct new feature
is twofold: the solutions concentrating around local minimum points of the potential $V$;
the nonlinearities are assumed to be either superlinear or asymptotically linear
at the infinity. Since the potential $V$ merely satisfies a local condition, we shall use the penalization technique to gain appropriate estimates of the solutions which ensure that the solutions are concentrated at local minimum points of the potential $V$. However, the penalized nonlinearity or the asymptotically linear nonlinearity results in some new difficulties as follows:
(i) for any $z\in E\backslash E^-$, where $E=H^1(\mathbb{R}^2)\times H^1(\mathbb{R}^2)$ and $E^-=\{(u,-u):u\in H^1(\mathbb{R}^2)\}$, the functional restricted in $\hat{E}(z)=\mathbb{R}z\oplus E^-$ may not have maximum points and then the method of generalized Nehari manifold which plays a vital role in \cite{ZhangCPDE} is invalid. Instead we make use of the reduction method, see \cite{AC,Ramos1} for example, and transform a strongly indefinite functional into a functional on $E^+$, where $E^+=\{(u,u): u\in H^1(\mathbb{R}^2)\}$. Unfortunately, the reduction functional has no maximum points either along the direction of some nontrivial $u\in H^1(\mathbb{R}^2)$. We shall use some techniques as in \cite{DINGXU1} and the relation between the system and its limit system to show that, for small enough $\epsilon>0$, the reduction functional has mountain-pass geometry, and the mountain-pass value has mini-max characterization which is helpful for showing the concentration of solutions as $\epsilon\rightarrow0$. (ii) It is not easy to prove that a family of solutions converge strongly  as $\epsilon\rightarrow0$, since there is no AR condition which plays an important role in \cite{ZhangCPDE}, we shall take advantage of a global compactness lemma, introduce another autonomous system and combine with some tricks to restore the compactness of solutions. (iii) Different from \cite{ZhangCPDE}, the asymptotically linear nonlinearity leads that the functional no longer satisfies (PS)$_c$ condition  and we shall show that in this case, the (Ce)$_c$ condition holds. In addition, the asymptotically linear nonlinearity enforces the implementation of some new estimates and verifications.

The paper is organized as
follows. In Section 2 we introduce the variational framework. In Section 3 we study two autonomous systems. In Section 4, we prove Theorems 1.1 and 1.2.

\section{Variational setting}
\renewcommand{\theequation}{2.\arabic{equation}}
In this paper we use the following
notations. Denote the norm in $L^r(\mathbb{R}^2)$
($1\leq r\leq\infty$) by $|\cdot|_r$. For simplicity, denote $\int_{\mathbb{R}^2}f(x)dx$ by $\int_{\mathbb{R}^2}f(x)$, denote $\frac{\partial h(x,y)}{\partial x}$ and $\frac{\partial h(x,y)}{\partial y}$ by $h'_1(x,y)$ and $h'_2(x,y)$ respectively. Without loss of generality, we assume $0\in \Lambda$ in condition (V$_2$).
If there is no special description in  lemmas, they are assumed that the functions $f$ and $g$ satisfy (H$_3$) or (H$'_3$). For subset $\Omega\subset\mathbb{R}^2$, denote $\mathbb{R}^2\backslash{\Omega}$ by $\Omega^c$.

Changing variable by $x\rightarrow\epsilon x$, problem (\ref{1.1})
turns out to be
\begin{equation}\label{2.1}\aligned
\left\{ \begin{array}{lll}
-\Delta u+V(\epsilon x)u=g(v)\ & \text{in}\quad \mathbb{R}^2,\\
-\Delta v+V(\epsilon x)v=f(u)\ & \text{in}\quad \mathbb{R}^2.
\end{array}\right.\endaligned
\end{equation}
Let $H^1(\mathbb{R}^2)$ be the Sobolev space endowed with the inner product and norm
$$(u,v)_\epsilon=\int_{\mathbb{R}^2}(\nabla u\nabla v+V(\epsilon x)uv),\ \|u\|^2_{\epsilon}=(u,u)_{\epsilon}, \ u,v\in H^1(\mathbb{R}^2),$$
and $E:=H^1(\mathbb{R}^2)\times H^1(\mathbb{R}^2)$  be the Sobolev space endowed with the inner product
$$(z_1,z_2)_{1,\epsilon}=(u_1,u_2)_{\epsilon}+(v_1,v_2)_{\epsilon},\quad z_i=(u_i,v_i)\in E, i=1,2.$$
It is easy to see that there is a space decomposition of $E$ that $E=E^+\oplus E^-$, where
$$E^+=\{(u,u): u\in H^1(\mathbb{R}^2)\}, \ E^-=\{(u,-u): u\in H^1(\mathbb{R}^2)\}.$$
For each $z=(u,v)\in E$, one has
$$z=z^++z^-=\bigl((u+v)/2,(u+v)/2\bigr)+\bigl((u-v)/2,(v-u)/2\bigr).$$
The energy functional of system (\ref{2.1}) is
$$I_\epsilon(z)=\int_{\mathbb{R}^2}(\nabla u\nabla v+V(\epsilon x)uv)-\int_{\mathbb{R}^2}(F(u)+G(v)).$$
Observe that
\begin{equation}\label{2.2}
I_\epsilon(z)=\frac12\|z^+\|^2_{1,\epsilon}
-\frac12\|z^-\|^2_{1,\epsilon}-\int_{\mathbb{R}^2}(F(u)+G(v)),
\end{equation}
which emphasizes the strongly indefinite nature of $I_\epsilon$.
\subsection{The modified problem}
Due to the local condition (V$_2$) on $V$, we apply the penalization method introduced by del Pino and Felmer \cite{DEF} to problem (\ref{2.1}). Fix small numbers $a_1, a_2>0$ in such a way that
$f'(a_1)\leq\frac{\inf_{\mathbb{R}^2}V}{2}$, $f'(t)\geq f'(a_1)$ for any $t\geq a_1$, $g'(a_2)\leq\frac{\inf_{\mathbb{R}^2}V}{2}$ and $g'(t)\geq g'(a_2)$ for any $t\geq a_2$.
Set
$$\aligned\tilde{f}(t)=\left\{\begin{array}{lll}f(t),\ \ &\text{if} \ t\leq a_1;\\ f'(a_1)t+f(a_1)-f'(a_1)a_1,\ \ &\text{if} \ t>a_1,\end{array}\right.\endaligned$$
$$\aligned\tilde{g}(t)=\left\{\begin{array}{lll}g(t),\ \ &\text{if} \ t\leq a_2;\\ g'(a_2)t+g(a_2)-g'(a_2)a_2,\ \ &\text{if} \ t>a_2.\end{array}\right.\endaligned$$
Then we introduce
$$\bar{f}(x,t)=\chi_{\Lambda}(x)f(t)
+(1-\chi_{\Lambda}(x))\tilde{f}(t),\quad \bar{g}(x,t)=\chi_{\Lambda}(x)g(t)
+(1-\chi_{\Lambda}(x))\tilde{g}(t),$$
and $\bar{F}(x,t)=\int^t_0\bar{f}(x,s)ds$, $\bar{G}(x,t)=\int^t_0\bar{g}(x,s)ds$. The relevant properties of $\bar{f}$ and $\bar{g}$ are
displayed in the next lemma, whose proof is elementary.
\begin{lemma}\label{l2.1} The function $\bar{f}(x,t)$ ( also $\bar{g}(x,t)$) satisfies:\\
\noindent(H$''_1$) $\bar{f}(x,t)=o(t)$\ uniformly in\
$x\in\mathbb{R}^2$, and  $\bar{f}(x,t)\leq f(t)$ for all $x\in\mathbb{R}^2$ and $t>0$;\\
\noindent(H$''_2$) (i) if (H$_3$) is satisfied, then $0< \theta \bar{F}(x,t)\leq \bar{f}(x,t)t,$ for all $x\in \Lambda$ and $t>0$;\\
\  \ \indent\  (ii) if (H$'_3$) is satisfied, then
 $0< 2 \bar{F}(x,t)< \bar{f}(x,t)t,$ for all $x\in \Lambda$ and $t>0$;\\
\noindent(H$''_3$) $0< 2\bar{F}(x,t)<\bar{f}(x,t)t\leq \frac{\inf_{\mathbb{R}^2}V}{2}t^2,$ for all $x\not\in\Lambda$ and $t>0$;\\
\noindent(H$''_4$) $\bar{f}(x,t)$ are nondecreasing in $t\in (0,+\infty)$ for all $x\in\mathbb{R}^2$; \\
\noindent(H$''_5$) set $$\hat{F}(x,t)=\frac12\bar{f}(x,t)t-\bar{F}(x,t),\quad \forall (x,t)\in\mathbb{R}^2\times(0,+\infty).$$ Then $\hat{F}(x,t)$ is nondecreasing in $t\in (0,+\infty)$ for all $x\in\mathbb{R}^2$ and  $\hat{F}(x,t)\rightarrow+\infty$ uniformly in $x$ as $t\rightarrow+\infty$;\\
\noindent(H$''_6$) for some (arbitrarily small) $\delta=\delta(a_1,a_2)>0$
\begin{equation}\label{2.2.1}|\bar{f}'_2(x,t)|=|\tilde{f}'(t)|\leq \delta, \quad\text{for all}\ x\not\in{\Lambda}\ \text{and}\ t>0.\end{equation}
\end{lemma}

Now we establish the modified problem
\begin{equation}\label{2.3}\aligned
\left\{ \begin{array}{lll}
-\Delta u+V(\epsilon x)u=\bar{g}(\epsilon x, v)\ & \text{in}\quad \mathbb{R}^2,\\
-\Delta v+V(\epsilon x)v=\bar{f}(\epsilon x, u)\ & \text{in}\quad \mathbb{R}^2.
\end{array}\right.\endaligned
\end{equation}
Denote $\Lambda_\epsilon=\{x\in\mathbb{R}^2:\epsilon x\in \Lambda\}$, then the solution $(u,v)$ of (\ref{2.3}) with $u(x)\leq a_1$ and $v(x)\leq a_2$ for each $x\in {\mathbb{R}^2\backslash{\Lambda_\epsilon}}$ is also a solution of (\ref{2.1}). The functional of (\ref{2.3}) is
$${\Phi}_\epsilon(z)=\int_{\mathbb{R}^2}(\nabla u\nabla v+V(\epsilon x)uv)-\int_{\mathbb{R}^2}\bar{G}(\epsilon x,v)-\int_{\mathbb{R}^2}\bar{F}(\epsilon x,u),\quad \forall z=(u,v)\in E,$$
 and ${\Phi}'_\epsilon\in C^1(E, \mathbb{R})$.

\begin{lemma}\label{l2.2} Let (V$_1$) and (H$_1$)-(H$_3$) hold. Then for  each $\epsilon>0$, the (PS)$_c$ condition of $\Phi_\epsilon$ with $c\in\mathbb{R}$ holds. \end{lemma}
{\bf Proof}: Let $\{z_n=(u_n,v_n)\}$ be the (PS)$_c$ sequence of $\Phi_\epsilon$. Firstly we show that $\{z_n\}$ is bounded in $E$. Note that \begin{equation}\label{3.14}\aligned c+o_n(1)\|z_n\|_{1,\epsilon}&=\Phi_\epsilon(z_n)-\frac12\langle \Phi'_\epsilon(z_n),z_n\rangle=\int_{\mathbb{R}^2}\hat{F}(\epsilon x,u_n)+\int_{\mathbb{R}^2}\hat{G}(\epsilon x,v_n)\\&\geq\bigl(\frac12-\frac1\theta\bigr)
\int_{\Lambda_\epsilon}\bigl(f(u_n)u_n+g(v_n)v_n\bigr).
\endaligned\end{equation}
By $\langle \Phi'_\epsilon(z_n),(u_n,0)\rangle=o_n(1)\|u_n\|_\epsilon$
and $\langle \Phi'_\epsilon(z_n),(0,v_n)\rangle=o_n(1)\|v_n\|_\epsilon$ we have
\begin{equation}\label{3.1.1} \|u_n\|_\epsilon=\int_{\mathbb{R}^2}\bar{g}(\epsilon x,v_n)\frac{u_n}{\|u_n\|_\epsilon}+o_n(1),\quad
\|v_n\|_\epsilon=\int_{\mathbb{R}^2}\bar{f}(\epsilon x,u_n)\frac{v_n}{\|v_n\|_\epsilon}+o_n(1).\end{equation}
In view of (H$''_1$), there exist $\beta>0$ and $C_\beta>0$ such that
\begin{equation*}
|\bar{f}(\epsilon x,t)|,|\bar{g}(\epsilon x,t)|\leq C_\beta e^{\beta t^2}\ \text{for all}\ (x,t)\in \mathbb{R}^2\times\mathbb{R}.
\end{equation*}
Set
$$\Lambda^1_{n,\epsilon}=\{x\in \Lambda_\epsilon: |\bar{f}(\epsilon x, u_n)|/ C_\beta\geq e^{\frac14}\},\quad \Lambda^2_{n,\epsilon}=\{x\in \Lambda_\epsilon: |\bar{f}(\epsilon x, u_n)|/ C_\beta\leq e^{\frac14}\}.$$
By (H$''_1$) there exists $C_1>0$ such that for all $n$,
$|\bar{f}(\epsilon x, u_n)|\leq C_1|u_n|$, for $x\in \Lambda^2_{n,\epsilon}$. As in \cite[Lemma 3.2]{DDR}, the following inequality holds
\begin{equation*}\aligned
st\leq
\left\{ \begin{array}{lll}
(e^{t^2}-1)+|s|(\log|s|)^{\frac12},\ & \ t\in\mathbb{R}\text{ and }\ |s|\geq e^{\frac14};\\
(e^{t^2}-1)+\frac12s^2,\ & \ t\in\mathbb{R}\text{ and }\ |s|\leq e^{\frac14}.
\end{array}\right.\endaligned
\end{equation*}
Applying the above inequality with $t=\frac{v_n}{\|v_n\|_\epsilon}$ and $s=\frac{\bar{f}(\epsilon x, u_n)}{C_\beta}$, from Trudinger-Moser inequality we get
$$\aligned\Bigl|\int_{\Lambda_\epsilon}\bar{f}(\epsilon x,u_n)\frac{v_n}{\|v_n\|_\epsilon}\Bigr|
\leq& C_\beta\Bigg[\int_{\Lambda^1_{n,\epsilon}}\frac{|\bar{f}(\epsilon x,u_n)|}{C_\beta}\Bigl[\log\bigl(\frac{|\bar{f}(\epsilon x,u_n)|}{C_\beta}\bigr)\Bigr]^{\frac12}
+\frac12\int_{\Lambda^2_{n,\epsilon}}\frac{1}{C^2_\beta}\bar{f}^2(\epsilon x,u_n)\Bigg]\\
&\hspace{4mm}+C_\beta
\int_{\Lambda_\epsilon}\bigl[e^{\bigl(\frac{v_n}{\|v_n\|}\bigr)^2}-1
\bigr]
\\ \leq& C_2+\bigl(\beta^{\frac12}+\frac{C_1}{2C_\beta}\bigr)
\int_{\Lambda_\epsilon}f(u_n)u_n.\endaligned$$
Moreover,
$$\Bigl|\int_{\Lambda^c_\epsilon}\bar{f}(\epsilon x,u_n)\frac{v_n}{\|v_n\|_\epsilon}\Bigr|
\leq\int_{\Lambda^c_\epsilon}\frac{\inf_{\mathbb{R}^2}V}{2}\frac{u_n v_n}{\|v_n\|_\epsilon}\leq\frac{\|u_n\|_\epsilon}{2}.$$
From (\ref{3.1.1}) it follows that
$$\|v_n\|_\epsilon\leq C_2+\bigl(\beta^{\frac12}+\frac{C_1}{2C_\beta}\bigr)\int_{\Lambda_\epsilon}f(u_n)u_n
+\frac{1}{2}\|u_n\|_\epsilon.$$
Similarly
$$\|u_n\|_\epsilon\leq C_2+\bigl(\beta^{\frac12}+\frac{C_1}{2C_\beta}\bigr)\int_{\Lambda_\epsilon}
g(v_n)v_n
+\frac{1}{2}\|v_n\|_\epsilon.$$
Using (\ref{3.14}) we get
$$\frac{1}{2}(\|u_n\|_\epsilon+\|v_n\|_\epsilon)\leq C_2+C\bigl(\beta^{\frac12}+\frac{C_1}{2C_\beta}\bigr)
\bigl(c+o_n(1)\|z_n\|_{1,\epsilon}\bigr).$$
Therefore, $\{z_n\}$ is bounded in $E$ and we assume that $z_n\rightharpoonup z_0=(u_0,v_0)$ in $E$.

Next we show that $z_n\rightarrow z_0$ in $E$. Since $g'$ and $f'$ are of subcritical exponential growth, it is easy to see that $g$ and $f$ are also of subcritical exponential growth. Then by (H$''_1$), for any fixed $\alpha,\beta>0$, then for any $\delta>0$ and any $q\geq1$, there exist $C_1(\delta,q)$ and  $C_2(\delta,q)$ such that
\begin{equation}\label{2.1.0}\aligned |\bar{f}(\epsilon x,t)|&\leq |f(t)|\leq \delta|t|+C_1(\delta,q)|t|^{q-1}(e^{\alpha t^2}-1),\quad \forall(x,t)\in \mathbb{R}^2\times\mathbb{R},\\
 |\bar{g}(\epsilon x,t)|&\leq|g(t)|\leq\delta|t|+C_2(\delta,q)|t|^{q-1}(e^{\beta t^2}-1), \quad \forall(x,t)\in \mathbb{R}^2\times\mathbb{R}.\endaligned\end{equation}
 In view of \cite[Lemma 2.1]{DDR} and (\ref{2.1.0}), one easily has that  $\Phi'_\epsilon$ is  weakly
sequentially continuous. Then $\Phi'_\epsilon(z_0)=0$. Since $\Lambda_\epsilon$ is bounded for any fixed $\epsilon>0$, using (\ref{2.1.0}) one easily has that
$$\int_{\Lambda_\epsilon}\bigl({f}(u_n )-{f}(u_0 )\bigr)(v_n-v_0)=o_n(1),\quad\int_{\Lambda_\epsilon}\bigl({g}(v_n )-{g}(v_0 )\bigr)(u_n-u_0)=o_n(1).$$
In addition, by (H$''_6$) we infer
$$\aligned&\Bigl|\int_{\Lambda^c_\epsilon}\bigl(\bar{g}(\epsilon x,v_n )-\bar{g}(\epsilon x,v_0 )\bigr)(u_n-u_0)\Bigr|=\Bigl|\int_{\Lambda^c_\epsilon}\bar{g}'_2(\epsilon x, \xi_n)(v_n-v_0)(u_n-u_0)\Bigr|\\ \leq& \int_{\Lambda^c_\epsilon} \frac{\inf_{\mathbb{R}^2}V}{2}|v_n-v_0||u_n-u_0|\leq \frac{\inf_{\mathbb{R}^2}V}{4}(|v_n-v_0|^2_2+|u_n-u_0|^2_2).\endaligned$$
Similarly for $\bar{f}$. Then
$$\aligned o_n(1)=&\langle \Phi'_\epsilon(z_n)-\Phi'_\epsilon(z_0),z_n-z_0\rangle\\
\geq&\|u_n-u_0\|^2_\epsilon+\|v_n-v_0\|^2_\epsilon-
\frac{\inf_{\mathbb{R}^2}V}{2}(|v_n-v_0|^2_2+|u_n-u_0|^2_2).\endaligned$$
Thus $z_n\rightarrow z_0$ in $E$.
\ \ \ \ \ $\Box$

\begin{lemma}\label{l6.1} Let (V$_1$), (H$_1$), (H$_2$) and (H$'_3$) hold. Then for each $\epsilon>0$, the (Ce)$_c$ condition of $\Phi_\epsilon$ with $c\in\mathbb{R}$ holds.
\end{lemma}
{\bf Proof}:
  Let $\{z_n=(u_n,v_n)\}$ be the (Ce)$_c$ sequence of $\Phi_\epsilon$ with $c\in\mathbb{R}$. Firstly, we show that $\{z_n\}$ is bounded in $E$. Argue by contradiction we assume $\|z_n\|_{1,\epsilon}\rightarrow+\infty$. Observe that
\begin{equation}\label{6.1} c+o_n(1)=\Phi_\epsilon(z_n)-\frac12\langle \Phi'_\epsilon(z_n),z_n\rangle=\int_{\mathbb{R}^2}\hat{F}(\epsilon x,u_n)+\int_{\mathbb{R}^2}\hat{G}(\epsilon x,v_n).\end{equation}
To get a contradiction, we set
$$\aligned &d_1(r)=\inf\{\hat{F}(\epsilon x,s):x\in\mathbb{R}^2 \ \text{and}\ s>r\},\ d_2(r)=\inf\{\hat{G}(\epsilon x,s):x\in\mathbb{R}^2\ \text{and}\ s>r\},\\ &\Omega^1_n(\rho,r)=\{x\in\mathbb{R}^2:\rho\leq |u_n(x)|\leq r\},\ \Omega^2_n(\rho,r)=\{x\in\mathbb{R}^2:\rho\leq |v_n(x)|\leq r\}\endaligned$$
 and
 $$\aligned C^r_\rho&=\inf\bigl\{\frac{\hat{F}(\epsilon x,s)}{s^2}:x\in\mathbb{R}^2,\ \rho\leq s\leq r\bigr\},
\bar{C}^r_\rho=\inf\bigl\{\frac{\hat{G}(\epsilon x,s)}{s^2}:x\in\mathbb{R}^2,\ \rho\leq s\leq r\bigr\}.
 \endaligned$$
By (H$''_5$), $d_1(r)\rightarrow+\infty$ as $r\rightarrow+\infty$. In view of the definition of  $C^r_\rho$ and (\ref{6.1}) we have
$$C\geq\int_{\Omega^1_n(0,\rho)}\hat{F}(\epsilon x,u_n(x))+C^r_\rho
\int_{\Omega^1_n(\rho,r)}|u_n(x)|^2+d_1(r)|\Omega^1_n(r,\infty)|.
$$
Then $|\Omega^1_n(r,\infty)|\rightarrow0$ as $r\rightarrow\infty$ uniformly in $n$, and
$\int_{\Omega^1_n(\rho,r)}|u_n(x)|^2\leq \frac{C}{C^r_\rho}.$ Similarly, $|\Omega^2_n(r,\infty)|\rightarrow0$ as $r\rightarrow\infty$ uniformly in $n$, $\int_{\Omega^2_n(\rho,r)}|v_n(x)|^2\leq \frac{C}{\bar{C}^r_\rho}.$
For any fixed $\delta>0$, by (H$''_1$) we know
$$|\bar{f}(\epsilon x,s)|\leq \delta|s|\ \text{and}\ |\bar{g}(\epsilon x,s)|\leq \delta|s|, \  \forall\ s\in[0,\rho_\delta].$$
Then
\begin{equation}\label{6.2} \Bigl|\int_{\Omega^1_n(0,\rho_{\delta})}\frac{\bar{f}(\epsilon x,u_n)v_n}{\|v_n\|^2_\epsilon}\Bigr|\leq C\delta\frac{\|u_n\|_\epsilon}{\|v_n\|_\epsilon},\ \ \Bigl|\int_{\Omega^2_n(0,\rho_\delta)}\frac{\bar{g}(\epsilon x,v_n)u_n}{\|u_n\|^2_\epsilon}\Bigr|\leq C\delta\frac{\|v_n\|_\epsilon}{\|u_n\|_\epsilon}.\end{equation}
To go on, up to a subsequence we may assume that $\frac{\|u_n\|_\epsilon}{\|v_n\|_\epsilon}\leq1$, since we can deal with the second inequality of (\ref{6.2}) to get a contradiction if $\frac{\|u_n\|_\epsilon}{\|v_n\|_\epsilon}>1$. By (\ref{6.2}) and $\frac{\|u_n\|_\epsilon}{\|v_n\|_\epsilon}\leq1$ we obtain
$\Bigl|\int_{\Omega^1_n(0,\rho_\delta)}\frac{\bar{f}(\epsilon x,u_n)v_n}{\|v_n\|^2_\epsilon}\Bigr|\leq C\delta.$
Since $\|z_n\|_{1,\epsilon}\rightarrow+\infty$ and $\frac{\|u_n\|_\epsilon}{\|v_n\|_\epsilon}\leq1$, we get $\|v_n\|_\epsilon\rightarrow+\infty$. Then
\begin{equation*}
\Bigl|\int_{\Omega^1_n(\rho_\delta,r_\delta)}\frac{\bar{f}(\epsilon x,u_n)v_n}{\|v_n\|^2_\epsilon}\Bigr|\leq
\frac{\int_{\Omega^1_n(\rho_\delta,r_\delta)}
l_0|u_n||v_n|}{\|v_n\|^2_\epsilon}\leq \frac{C}{C^{r_\delta}_{\rho_\delta}\|v_n\|_\epsilon}\rightarrow0.
\end{equation*}
Moreover,
\begin{equation*}
\Bigl|\int_{\Omega^1_n(r_\delta,+\infty)}\frac{\bar{f}(\epsilon x,u_n)v_n}{\|v_n\|^2_\epsilon}\Bigr|\leq
\frac{l_0|u_n|_3|v_n|_3|\Omega^1_n(r_\delta,\infty)|^{\frac13}}
{\|v_n\|^2_\epsilon}\leq C\delta.
\end{equation*}
Then $\int_{\mathbb{R}^2}\frac{\bar{f}(\epsilon x,u_n)v_n}{\|v_n\|^2_\epsilon}\rightarrow0$.
On the other hand, note that
$$o_n(1)=\langle \Phi'_\epsilon(u_n,v_n),(0,v_n)\rangle=\|v_n
\|^2_\epsilon-
\int_{\mathbb{R}^2}\bar{f}(\epsilon x,u_n)v_n,$$
and $\|v_n\|_\epsilon\rightarrow+\infty$, we infer
$\int_{\mathbb{R}^2}\frac{\bar{f}(\epsilon x,u_n)v_n
}{\|v_n\|^2_\epsilon}\rightarrow1.$
This is a contradiction.
Therefore, $\|z_n\|_{1,\epsilon}$ is bounded. Taking same arguments as in Lemma \ref{l2.2} we know $\{z_n\}$ converges strongly in $E$.\ \ \ \ $\Box$

 \subsection{Functional reduction}
Since the functional $\Phi_\epsilon$ is strongly indefinite, we shall apply the reduction approaches, see for example \cite{AC,Ramos1}, to look for critical  points of $\Phi_\epsilon$. More precisely, we shall  reduce
the strongly indefinite functional $\Phi_\epsilon$ to a
functional on $E^+$.

 For any fixed $(u,u)\in E^+$, let $\phi_{(u,u)}:E^-\rightarrow\mathbb{R}$ defined by $$\phi_{(u,u)}(v,-v)=\Phi_\epsilon(u+v,u-v).$$ Then
\begin{equation}\label{3.1}\aligned
\phi_{(u,u)}(v,-v)
&=\|u\|^2_\epsilon-\|v\|^2_\epsilon-\int_{\mathbb{R}^2}[\bar{F}(\epsilon x,u+v)+\bar{G}(\epsilon x,u-v)]\leq\|u\|^2_\epsilon-\|v\|^2_\epsilon.
\endaligned
\end{equation}
Moreover, for any $(v,-v)$, $(w,-w)\in E^-$, by (H$''_4$) we have
\begin{equation}\label{3.2}\aligned
\phi''_{(u,u)}[(v,-v),(w,-w)]&=-2\|w\|^2_\epsilon-
\int_{\mathbb{R}^2}[\bar{g}'_2(\epsilon x,u-v)+\bar{f}'_2(\epsilon x,u+v)]w^2\leq-2\|w\|^2_\epsilon.
\endaligned
\end{equation}
Due to (\ref{3.1}) and (\ref{3.2}), there exists a unique $\bar{h}_\epsilon(u)\in H^1(\mathbb{R}^2)$ such that
$$\Phi_\epsilon(u+\bar{h}_\epsilon(u),u-\bar{h}_\epsilon(u))=\max_{v\in H^1(\mathbb{R}^2)}\Phi_\epsilon(u+v,u-v).$$
Consequently, the operator $\bar{h}_\epsilon: H^1(\mathbb{R}^2)\rightarrow H^1(\mathbb{R}^2)$ is well defined, and
\begin{equation}\label{3.4}\langle\Phi'_\epsilon(u+\bar{h}_\epsilon(u),
u-\bar{h}_\epsilon(u)),(\varphi,-\varphi)\rangle=0, \quad \forall \varphi\in H^1(\mathbb{R}^2).\end{equation}

\begin{lemma}\label{l2.4} The map $\bar{h}_\epsilon$ belongs to $C^1(H^1(\mathbb{R}^2), H^1(\mathbb{R}^2))$.
\end{lemma}
{\bf Proof}: Define the map
$$\mathcal{H}:E\times E^-\rightarrow E^-,\ \mathcal{H}((u,v),(\psi,-\psi))=P\circ\Phi'_\epsilon(u+\psi,u-\psi),$$
where $P$ is the projection from $E$ to $E^-$.
Note that $\mathcal{H}$ is class $C^1$ and its partial derivative of $\mathcal{H}$ with respect to the second variable is
$$\mathcal{H}'_2((u,v),(\psi,-\psi))(\varphi,-\varphi)=P\circ \Phi''_\epsilon(u+\psi,v-\psi)(\varphi,-\varphi).$$
If we identify by $E^-=(E^-)^*$ and define $T=P\circ \Phi''_\epsilon(u+\psi,v-\psi),$
then
$$\langle T(\phi,-\phi),(\varphi,-\varphi)\rangle=-2(\phi,\varphi )_\epsilon-\int_{\mathbb{R}^2}[\bar{f}'_2(\epsilon x,u+\psi)+\bar{g}'_2(\epsilon x,u-\psi)]\phi\varphi,\quad\forall (\phi,\varphi)\in E.$$
We claim that $T$ is one-to-one. In fact, if $T(\phi,-\phi)=(0,0)$, then
$\langle T(\phi,-\phi),(\phi,-\phi) \rangle=0$. Using (H$''_4$) we have $\phi=0$. On the other hand, we have
$$\aligned (id+T)((\phi,-\phi),(\varphi,-\varphi))&=
\int_{\mathbb{R}^2}\bar{f}'_2(\epsilon x,u+\psi)\phi\varphi+\int_{\mathbb{R}^2}\bar{g}'_2(\epsilon x,u-\psi)\phi\varphi\\&:=T_1((\phi,-\phi),(\varphi,-\varphi))+T_2((\phi,-\phi),(\varphi,-\varphi)),\endaligned$$
where $id: E^-\rightarrow(E^-)^*$, $id((\phi,-\phi),(\varphi,-\varphi))=2(\phi,\varphi)_\epsilon$ for all
$\phi,\varphi\in H^1(\mathbb{R}^2)$ and $T_1,T_2:E^-\rightarrow(E^-)^*$. To show $id+T$ is a compact operator, we prove $T_1,T_2$ are compact.
In view of (H$''_1$) and (H$_3$)-(i), for any fixed $\alpha>0$, then for any $\delta>0$, there exists $C_\delta>0$ such that
\begin{equation}\label{2.2.0}|\bar{f}'_2(x,t)|\leq \delta+C_\delta(e^{\alpha t^2}-1), \ \forall (x,t)\in\mathbb{R}^2\times\mathbb{R}.\end{equation}
If $(\phi_n,-\phi_n)\rightharpoonup (\phi_0,-\phi_0)$ in $E^-$, then for some large $R>0$ we have
$$\aligned&\int_{\mathbb{R}^2}|\bar{f}'_2(\epsilon x,u+\psi)||\phi_n-\phi_0||\varphi|\leq \int_{\mathbb{R}^2}(\delta+C_\delta [e^{\alpha|u+\psi|^2}-1])|\phi_n-\phi_0||\varphi|\\
\leq&\delta|\phi_n-\phi_0|_2|\varphi|_2+C_\delta|\varphi|_2
\Bigl[\int_{|x|<R}\bigl(e^{4\alpha|u+\psi|^2}-1\bigr)\Bigr]^\frac14
\Bigl(\int_{|x|<R}|\phi_n-\phi_0|^4\Bigr)^{\frac14}\\
&+C_\delta|\varphi|_2\Bigl[\int_{|x|>R}\bigl(e^{4\alpha|u+\psi|^2}-1\bigr)\Bigr]^\frac14
\Bigl(\int_{|x|>R}|\phi_n-\phi_0|^4\Bigr)^{\frac14}
=o_n(1)\|\varphi\|_\epsilon.\endaligned$$
Hence $$\|T_1(\phi_n, -\phi_n)-T_1(\phi_0, -\phi_0)\|_{(E^-)^*}\rightarrow0.$$
Then $T_1$ is compact. Similarly $T_2$ is compact and so $id+T$ is a compact operator. From the Fredholm alternative we conclude that $T$ is isomorphism. From (\ref{3.4}) we know
$$P\circ \Phi'_\epsilon(w+\bar{h}_\epsilon(w), w-\bar{h}_\epsilon(w))=0,\quad\ \forall w\in H^1(\mathbb{R}^2).$$
Thus the implicit function theorem implies that $\bar{h}_\epsilon$ is of class $C^1$.\ \ \ \ $\Box$

By Lemma \ref{l2.4}, we consider the reduced functional $J_\epsilon: H^1(\mathbb{R}^2)\rightarrow\mathbb{R}$ defined by
$$\aligned J_\epsilon(w)&=
\Phi_\epsilon(w+\bar{h}_\epsilon(w),w-\bar{h}_\epsilon(w))\\
&=\|w\|^2_\epsilon-\|\bar{h}_\epsilon(w)\|^2_\epsilon-\int_{\mathbb{R}^2}\bar{F}(\epsilon x, w+\bar{h}_\epsilon(w))-\int_{\mathbb{R}^2}\bar{G}(\epsilon x, w-\bar{h}_\epsilon(w)),\endaligned$$
which is of class $C^1$. Moreover, by (\ref{3.4}) we have
\begin{equation}\label{2.2.2}\aligned\langle J'_\epsilon(w),\phi\rangle&= \langle\Phi'_\epsilon(w+\bar{h}_\epsilon(w),w-\bar{h}_\epsilon(w)),
(\phi+\bar{h}'_\epsilon(w)\phi,\phi-\bar{h}'_\epsilon(w)\phi)\rangle\\
&= \langle\Phi'_\epsilon(w+\bar{h}_\epsilon(w),w-\bar{h}_\epsilon(w)),
(\phi,\phi)\rangle, \quad \forall\phi\in H^1(\mathbb{R}^2).\endaligned\end{equation}

\begin{lemma}\label{l2.5}(1) The map $$\eta_\epsilon: H^1(\mathbb{R}^2)\rightarrow E: u\rightarrow (u+\bar{h}_\epsilon(u), u-\bar{h}_\epsilon(u))$$
is a homeomorphism between critical points of $J_\epsilon$ and $\Phi_\epsilon$ and $\eta^{-1}_\epsilon: E\rightarrow H^1(\mathbb{R}^2)$ is given by $\eta^{-1}_\epsilon(u,v)=\frac{u+v}{2}$.\\
\noindent(2) If $\{w_n\}$ is a (PS)$_c$ sequence with any $c\in\mathbb{R}$ for $J_\epsilon$, then $\{(w_n+\bar{h}_\epsilon(w_n),w_n-\bar{h}_\epsilon(w_n))\}$ is
a (PS)$_c$ sequence for $\Phi_\epsilon$. \\
\noindent(3) If $f$ and $g$ satisfy (H$'_3$), and  $\{w_n\}$ is a (Ce)$_c$ sequence with any $c\in\mathbb{R}$ for $J_\epsilon$, then $\{(w_n+\bar{h}_\epsilon(w_n),w_n-\bar{h}_\epsilon(w_n))\}$ is
a (Ce)$_c$ sequence for $\Phi_\epsilon$.
\end{lemma}
{\bf Proof}: The proof of (1) is easy and we omit. For any $(\phi,\varphi)\in E$, by (\ref{3.4}) we get
$$\aligned&\langle \Phi'_\epsilon(w_n+\bar{h}_\epsilon(w_n),w_n-\bar{h}_\epsilon(w_n)),
(\phi,\varphi)\rangle
\\=&\langle \Phi'_\epsilon(w_n+\bar{h}_\epsilon(w_n),w_n-\bar{h}_\epsilon(w_n)),
(\frac{\phi+\varphi}{2}, \frac{\phi+\varphi}{2})\rangle=\langle J'_\epsilon(w_n),\frac{\phi+\varphi}{2}\rangle.
\endaligned$$
Then the conclusion (2) holds true. Regarding (3), firstly  we introduce an inequality
\begin{equation}\label{3.1.3}\aligned 0\leq& \Phi_\epsilon(u+\bar{h}_\epsilon(u),u-\bar{h}_\epsilon(u))-\Phi_\epsilon(u,u)\\
=&\|u\|^2_\epsilon-\|\bar{h}_\epsilon(u)\|^2_\epsilon-
\int_{\mathbb{R}^2}\bar{F}(\epsilon x,u+\bar{h}_\epsilon(u))-\int_{\mathbb{R}^2}\bar{G}(\epsilon x,u-\bar{h}_\epsilon(u))\\&-\Bigl[\|u\|^2_\epsilon-\int_{\mathbb{R}^2}\bar{F}(\epsilon x,u)-\int_{\mathbb{R}^2}\bar{G}(\epsilon x,u)\Bigr]\\
\leq&-\|\bar{h}_\epsilon(u)\|^2_\epsilon+
\int_{\mathbb{R}^2}\bar{F}(\epsilon x,u)+\int_{\mathbb{R}^2}\bar{G}(\epsilon x,u).
\endaligned\end{equation}
 If (H$'_3$) is satisfied, i.e. $f,g$ are asymptotically linear at infinity, by (\ref{3.1.3}) we know
$$\|\bar{h}_\epsilon(u)\|^2_\epsilon\leq \int_{\mathbb{R}^2}\bar{F}(\epsilon x,u)+\int_{\mathbb{R}^2}\bar{G}(\epsilon x,u)\leq C\|u\|^2_\epsilon.$$
Then
$$\|(u+\bar{h}_\epsilon(u),u-\bar{h}_\epsilon(u))\|^2_{1,\epsilon}
=\|u+\bar{h}_\epsilon(u)\|^2_\epsilon+\|u-\bar{h}_\epsilon(u)\|^2_\epsilon\leq
C\|u\|^2_\epsilon,$$
Combining with the above (2),  the conclusion (3) yields. \ \ \  \ $\Box$

According to Lemma \ref{l2.5}, it suffices to look for critical points of $J_\epsilon$ and we shall show that $J_\epsilon$ possesses the mountain-pass structure.
\begin{lemma}\label{l2.6} There are $r>0$ and $\tau>0$ both independent of $\epsilon$, such that $J_\epsilon|_{S_r}\geq\tau$, where
$S_r=\{u\in H^1(\mathbb{R}^2):\|u\|_\epsilon=r\}$.\end{lemma}
{\bf Proof}: For any $w\in H^1(\mathbb{R}^2)$, notice that
$$J_\epsilon(w)=\Phi_\epsilon(w+\bar{h}_\epsilon(w), w-\bar{h}_\epsilon(w))=\max_{v\in H^1(\mathbb{R}^2)}\Phi_\epsilon(w+v, w-v)\geq\Phi_\epsilon(w, w).$$
By (H$''_2$), (\ref{2.1.0}), the H\"{o}lder inequality and Trudinger-Moser inequality we obtain
\begin{equation*}\aligned J_\epsilon(w)&\geq\|w\|^2_\epsilon-\int_{\mathbb{R}^2}\bar{F}(\epsilon x,w)-\int_{\mathbb{R}^2}\bar{G}(\epsilon x,w)\\
&\geq\|w\|^2_\epsilon-2\delta|w|^2_2-C_1(\delta,q)
|w|^q_{2q}\int_{\mathbb{R}^2}\bigl(e^{2\alpha w^2}-1\bigr)-C_2(\delta,q)
|w|^q_{2q}\int_{\mathbb{R}^2}\bigl(e^{2\beta w^2}-1\bigr)
\\&\geq\|w\|^2_\epsilon-2\delta|w|^2_2-C|w|^q_{2q}.
\endaligned\end{equation*}
Choosing $q>2$ in (\ref{2.1.0}), for some $r,\tau>0$ we have $J_\epsilon(w)\geq\tau>0$ when $\|w\|_\epsilon=r$. \ \ \ \ \ $\Box$

\begin{lemma}\label{l3.6} For any $u\in H^1(\mathbb{R}^2)\backslash\{0\}$, the following results hold.\\
\noindent(1) If (H$_3$) is satisfied, then\\
\ \indent(i) $J_\epsilon(tu)\rightarrow-\infty$ as $t\rightarrow+\infty$ if $suppu \cap \Lambda_\epsilon\neq\emptyset$;\\
\ \indent(ii) $J_\epsilon(tu)\rightarrow-\infty$ or $J_\epsilon(tu)\rightarrow+\infty$ as $t\rightarrow+\infty$ if $supp u\subset \mathbb{R}^2\backslash{\Lambda_\epsilon}$.\\
\noindent(2) If (H$'_3$) is satisfied, then $J_\epsilon(tu)\rightarrow-\infty$ or $J_\epsilon(tu)\rightarrow+\infty$ as $t\rightarrow+\infty$.
\end{lemma}

{\bf Proof}: (1) (i) Assume $supp u\cap \Lambda_\epsilon\neq\emptyset$.  Note that
$$\aligned J_\epsilon(tu)
&=t^2\|u\|^2_\epsilon-\|\bar{h}_\epsilon(tu)\|^2_\epsilon-\int_{\mathbb{R}^2}\bar{F}(\epsilon x,tu+\bar{h}_\epsilon(tu))-\int_{\mathbb{R}^2}\bar{G}(\epsilon x,tu-\bar{h}_\epsilon(tu))\\
&\leq t^2\|u\|^2_\epsilon-\|\bar{h}_\epsilon(tu)\|^2_\epsilon-
\int_{\Lambda_\epsilon}{F}(tu+\bar{h}_\epsilon(tu))-
\int_{\Lambda_\epsilon}{G}(tu-\bar{h}_\epsilon(tu)).\endaligned$$
By (H$''_2$)-(i) and (H$''_1$), for any $\delta>0$ there exist $C_{1,\delta}$ and $C_{2,\delta}>0$ such that
$$F(s)\geq C_{1,\delta} s^\theta-{\delta} s^2, \ G(s)\geq C_{2,\delta} s^\theta-{\delta} s^2,\ \forall s\in\mathbb{R}^+.$$
Then
$$\aligned J_\epsilon(tu)\leq& ({1+2\delta})t^2\|u\|^2_\epsilon-(1-2\delta)\|\bar{h}_\epsilon(tu)\|^2_\epsilon-
C_{1,\delta}\int_{\Lambda_\epsilon}|tu+\bar{h}_\epsilon(tu)|^\theta\\&-
C_{2,\delta}\int_{\Lambda_\epsilon}|tu-\bar{h}_\epsilon(tu)|^\theta\\
\leq&({1+2\delta})t^2\|u\|^2_\epsilon-\min\{C_{1,\delta},C_{2,\delta}\}
\int_{\Lambda_\epsilon}2^\theta t^\theta|u|^\theta.
\endaligned$$
Hence
$J_\epsilon(tu)\rightarrow -\infty$ as $t\rightarrow+\infty$.

(ii) Suppose $supp u\subset \mathbb{R}^2\backslash{\Lambda_\epsilon}$.
If $J_\epsilon(tu)\rightarrow+\infty$, we are done. Otherwise,
we may assume $\sup_{t\geq0}J_\epsilon(tu)=M<+\infty$. Then
$$\aligned\frac{d}{dt}J_\epsilon(tu)&
=\frac1t\langle J'_\epsilon(tu),tu\rangle
=\frac1t\langle \Phi'_\epsilon(tu+\bar{h}_\epsilon(tu),tu-\bar{h}_\epsilon(tu)),
(tu+\bar{h}_\epsilon(tu),tu-\bar{h}_\epsilon(tu))\rangle
\\
&=\frac{2J_\epsilon(tu)}{t}-
\frac{2}{t}\int_{\mathbb{R}^2}\hat{F}(\epsilon x, tu+\bar{h}_\epsilon(tu))-\frac{2}{t}\int_{\mathbb{R}^2}\hat{G}(\epsilon x, tu-\bar{h}_\epsilon(tu)).\endaligned$$
For $r>0$, by (H$''_5$) and the fact that $f(t)=0$ if $t\leq0$, we infer
\begin{equation*}\aligned
\int_{\mathbb{R}^2}\hat{F}(\epsilon x, tu+\bar{h}_\epsilon(tu))\geq&
\int_{\{\mathbb{R}^2:u+{\bar{h}_\epsilon(tu)}/{t}\geq r\}}\hat{F}(\epsilon x, tu+\bar{h}_\epsilon(tu))\\
\geq& \hat{F}(\epsilon x,rt)meas\{x\in\mathbb{R}^2:u+{\bar{h}_\epsilon(tu)}/{t}\geq r\}.
\endaligned\end{equation*}
Similarly,
$$\aligned
\int_{\mathbb{R}^2}\hat{G}(\epsilon x, tu-\bar{h}_\epsilon(tu))\geq&
 \hat{G}(\epsilon x,rt)meas\{x\in\mathbb{R}^2:u-{\bar{h}_\epsilon(tu)}/{t}\geq r\}.
\endaligned$$
From (\ref{3.1.3}) and $supp u\subset\mathbb{R}^2\backslash{\Lambda_\epsilon}$ we know
\begin{equation*}\aligned\|\bar{h}_\epsilon(tu)\|^2_\epsilon&\leq \int_{\mathbb{R}^2\backslash{\Lambda_\epsilon}}\bigl(\bar{F}(\epsilon x,tu)+\bar{G}(\epsilon x,tu)\bigr)\leq \frac{\inf_{\mathbb{R}^2}V}{2}t^2|u|^2_2.\endaligned\end{equation*}
Then $\{\frac{\bar{h}_\epsilon(tu)}{t}\}_{t>0}\subset H^1(\mathbb{R}^2)$ is bounded. Since $u\neq0$, it is easy to see that one of the following two cases occurs \begin{equation*}meas\{x\in\mathbb{R}^2:
u-{\bar{h}_\epsilon(tu)}/{t}\geq r\}\geq \delta \ \text{or}\ meas\{x\in\mathbb{R}^2:u+{\bar{h}_\epsilon(tu)}/{t}\geq r\}\geq \delta,\end{equation*} with some $\delta>0$ for all $t>0$ provided $r>0$ small. Without loss of generality, we assume $meas\{x\in\mathbb{R}^2:u+{\bar{h}_\epsilon(tu)}/{t}\geq r\}\geq \delta$. By (H$''_5$) we get
$$\frac{d {J_\epsilon(tu)}}{dt}\leq \frac{2J_\epsilon(tu)}{t}-\frac{2\hat{F}(\epsilon x,rt)}{t}{\delta}\leq\frac{2J_\epsilon(tu)}{t}-\frac{3M}{t}\leq-\frac{M}{t}.$$
So
$J_\epsilon(tu)=\int^t_0\frac{dJ_\epsilon(tu)}{dt}\rightarrow-\infty$, as\ $t\rightarrow+\infty.$

(2) If (H$'_3$) is satisfied, then taking similar arguments as in the proof of the conclusion (1)-(ii), we obtain desired results.
\ \ \ \ \ $\Box$

\begin{lemma}\label{l3.7}
If $u\in H^1(\mathbb{R}^2)\backslash\{0\}$ satisfies $\langle J'_\epsilon(u),u\rangle=0$, then $J''_\epsilon(u)[u,u]<0$.
\end{lemma}

{\bf Proof}:
Observe that
\begin{equation*}\aligned J''_\epsilon(u)[u,u]=&2\|u\|^2_\epsilon-
2\|\bar{h}'_\epsilon(u)u\|^2_\epsilon
-\int_{\mathbb{R}^2}\bar{g}'_2(\epsilon x,u-\bar{h}_\epsilon(u))(u-\bar{h}'_\epsilon(u)u)^2\\
&-\int_{\mathbb{R}^2}\bar{f}'_2(\epsilon x,u+\bar{h}_\epsilon(u))(u+\bar{h}'_\epsilon(u)u)^2.\endaligned
\end{equation*}
Setting $$z_1=u+\bar{h}_\epsilon(u),\ z_2=u-\bar{h}_\epsilon(u), w_1=\bar{h}'_\epsilon(u)u-\bar{h}_\epsilon(u),\
$$ we obtain \begin{equation}\label{3.1.7}\aligned
J''_\epsilon(u)[u,u]
=&2\|u\|^2_\epsilon-2\bigl[(\bar{h}_\epsilon(u),\bar{h}_\epsilon(u)+2w_1)_\epsilon
+\|w_1\|^2_\epsilon\bigr]\\&-
\int_{\mathbb{R}^2}\bar{f}'_2(\epsilon x,z_1)(z_1+w_1)^2-\int_{\mathbb{R}^2}\bar{g}'_2(\epsilon x,z_2)(z_2-w_1)^2.
\endaligned\end{equation}
Note that \begin{equation}\label{3.1.4}\aligned
0=&\langle J'_\epsilon(u),u\rangle+\langle \Phi'_\epsilon(u+\bar{h}_\epsilon(u),u-\bar{h}_\epsilon(u)),(v,-v)\rangle\\
=&2\|u\|^2_\epsilon-2( \bar{h}_\epsilon(u),\bar{h}'_\epsilon(u)u+v)_\epsilon-
\int_{\mathbb{R}^2}f(\epsilon x,u+\bar{h}_\epsilon(u))(u+\bar{h}'_\epsilon(u)u+v)\\&-
\int_{\mathbb{R}^2}g(\epsilon x,u-\bar{h}_\epsilon(u))(u-\bar{h}'_\epsilon(u)u-v),\ \forall v\in H^1(\mathbb{R}^2).
\endaligned\end{equation}
Letting $v=w_1$ in (\ref{3.1.4}) we have
\begin{equation}\label{3.1.8}0=2\|u\|^2_\epsilon-2( \bar{h}_\epsilon(u),\bar{h}_\epsilon(u)+2w_1)_\epsilon-
\int_{\mathbb{R}^2}\bar{f}(\epsilon x,z_1)(z_1+2w_1)-\int_{\mathbb{R}^2}\bar{g}(\epsilon x,z_2)(z_2-2w_1).\end{equation}
Then (\ref{3.1.7}) minus (\ref{3.1.8}) we deduce
$$\aligned J''_\epsilon(u)[u,u]
=&-2\|w_1\|^2_\epsilon-\int_{\mathbb{R}^2}\frac{\bar{f}(\epsilon x,z_1)}{z_1}w^2_1+\int_{\mathbb{R}^2}\bigl[\frac{\bar{f}(\epsilon x,z_1)}{z_1}-\bar{f}'_2(\epsilon x,z_1)\bigr](z_1+w_1)^2\\
&-\int_{\mathbb{R}^2}\frac{\bar{g}(\epsilon x,z_2)}{z_2}w^2_1+\int_{\mathbb{R}^2}\bigl[
\frac{\bar{g}(\epsilon x,z_2)}{z_2}-
\bar{g}'_2(\epsilon x,z_2)\bigr](z_2-w_1)^2.\endaligned$$
Since $z_1+z_2=2u\neq0$, we may assume that $z_1\neq0$. So $J''_\epsilon(u)[u,u]<0$.\ \ \ \  \ $\Box$

\section{The autonomous systems}
\renewcommand{\theequation}{3.\arabic{equation}}
\subsection{The limit system}
For any $\mu\in[\inf_{\mathbb{R}^2}V,|V|_\infty]$, consider the limit system
\begin{equation}\label{5.1}\aligned
\left\{ \begin{array}{lll}
-\Delta u+\mu u=g(v)\ & \text{in}\quad \mathbb{R}^2,\\
-\Delta v+\mu v=f(u)\ & \text{in}\quad \mathbb{R}^2.
\end{array}\right.\endaligned
\end{equation}
To study (\ref{5.1}), for any $\mu>0$, define another inner product and norm in $H^1(\mathbb{R}^2)$ by
 $$(u,v)_\mu=\int_{\mathbb{R}^2}(\nabla u\nabla v+\mu uv),\ \|u\|^2_{\mu}=(u,u)_{\mu}, \ u,v\in H^1(\mathbb{R}^2),$$
 and another inner product in $E$ as follows $$(z_1,z_2)_{1,\mu}=(u_1,u_2)_{\mu}+(v_1,v_2)_{\mu},\quad z_i=(u_i,v_i)\in E, i=1,2.$$
The functional  of (\ref{5.1}) is
$$\Phi_\mu(z)=\frac12(\|z^+\|^2_{1,\mu}-\|z^-\|^2_{1,\mu})
-\int_{\mathbb{R}^2}[G(v)+F(u)], \ \forall z=z^++z^-\in E=E^+\oplus E^-.$$
Define
\begin{equation}\label{3.1.0}c_\mu:=\inf_{\mathcal{K}_\mu}\Phi_\mu,\quad \text{where}\ \mathcal{K}_\mu=\bigl\{ z\in E\backslash\{(0,0)\}, \ \Phi'_\mu(z)=0\bigr\}.\end{equation}
We would like to point out that similar to the arguments in \cite{ZhangCPDE} for the superlinear case and in \cite{LY} for the asymptotically linear case, we can show that  (\ref{5.1})  has a ground state and weaken $f$ and $g$ to be merely continuous.
In this paper, to give another characterization of the least energy $c_\mu$ to restore the compactness of (\ref{2.3}), we assume $f, g\in C^1$ and introduce the mapping $\bar{h}_\mu$.
As in Section 2.2, for any $u\in H^1(\mathbb{R}^2)$, there exists a unique $\bar{h}_\mu(u)\in H^1(\mathbb{R}^2)$ such that
\begin{equation}\label{3.3} \Phi_\mu(u+\bar{h}_\mu(u),u-\bar{h}_\mu(u))=\max_{v\in H^1(\mathbb{R}^2)}\Phi_\mu(u+v,u-v),\end{equation}
and define
$$J_\mu: H^1(\mathbb{R}^2)\rightarrow\mathbb{R},\ J_\mu(u)=\Phi_\mu(u+\bar{h}_\mu(u),u-\bar{h}_\mu(u)).$$
There is a one-to-one correspondence between the critical points of $\Phi_\mu$ and $J_\mu$. Similar to (\ref{3.1.3}), there holds
\begin{equation}\label{4.0.4}
\|\bar{h}_\mu(u)\|^2_\mu\leq
\int_{\mathbb{R}^2}{F}(u)+\int_{\mathbb{R}^2}{G}(u).\end{equation}

{\bf 3.1.1. The super-linear case.} In the same way as the proof of Lemma \ref{l2.6} and Lemma \ref{l3.6} (1)-(i) respectively,
we have the following results.
\begin{lemma}\label{l4.7}
(1)\ There are $r, \tau>0$ such that $J_\mu|_{S_r}\geq\tau$, where
$S_r=\{u\in H^1(\mathbb{R}^2):\|u\|_\mu=r\}$.\\
(2) For any $u\in H^1(\mathbb{R}^2)\backslash\{0\}$, $J_\mu(tu)\rightarrow-\infty$ as $t\rightarrow+\infty$.
\end{lemma}
\begin{lemma}\label{l5.1}
Let (H$_1$)-($H_3$) hold. Then (\ref{5.1}) admits a ground state in $E$. Moreover,
$c_\mu=\bar{c}_\mu=\bar{\bar{{c}}}_\mu,$
where \begin{equation*}\label{3.18.0}\bar{c}_\mu=\inf_{v\in \Gamma_\mu}\max_{t\in[0,1]}J_\mu(\nu(t))\ \text{with}\ \Gamma_\mu=\{\nu\in C([0,1], H^1(\mathbb{R}^2)):\nu(0)=0, J_\mu(\nu(1))<0\},\end{equation*} and
$$\bar{\bar{{c}}}_\mu=\inf_{u\in H^1(\mathbb{R}^2)\backslash\{0\}}\max_{t>0}J_\mu(tu).$$
\end{lemma}
{\bf Proof}:  In view of Lemma \ref{l4.7}, the mountain-pass theorem implies that $J_\mu$ has a (PS)$_{\bar{c}_\mu}$ sequence $\{w_n\}$ in $H^1(\mathbb{R}^2)$.
 Then Lemma \ref{l2.5}-(2) implies $\{z_n=(w_n+\bar{h}_\mu(w_n),w_n-\bar{h}_\mu(w_n))\}$ is a (PS)$_{\bar{c}_\mu}$  sequences of $\Phi_\mu$.
Similar to the proof of Lemma \ref{l2.2}, $\{z_n\}$ is bounded in $E$ and assume $\|z_n\|^2_{1,\mu}\leq M$. Below we claim that $\{z_n\}$ is non-vanishing. Otherwise if $\{z_n\}$ is vanishing, i.e.
 $$\lim_{n\rightarrow\infty}\sup_{y\in\mathbb{R}^2}\int_{B_R(y)}
 (u^2_n+v^2_n)=0, \ \ \forall R>0,$$
then P.L. Lions compactness Lemma implies that $u_n\rightarrow0$ and $v_n\rightarrow0$ in $L^r(\mathbb{R}^2)$ for any $r>2$.
For some $\alpha\in (0,\frac{8\pi}{3M})$, by (\ref{2.1.0}) with $q=1$, for any $\delta>0$, there exists $C(\delta)>0$ such that
\begin{equation}\label{5.3}\aligned
\Bigl|\int_{\mathbb{R}^2}{f}(u_n)v_n\Bigr|&
&\leq\delta|u_n|_2|v_n|_2+C(\delta)
|v_n|^{3}_{3}\Bigl(\int_{\mathbb{R}^2}
\bigl[e^{\alpha \frac{3}{2} \|u_n\|^2_\mu\frac{u^2_n}{\|u_n\|^2_\mu}}-1\bigr]\Bigr)^{\frac{2}{3}}\rightarrow0.
\endaligned\end{equation}
Similarly, we get
$\int_{\mathbb{R}^2}{g}(v_n)u_n\rightarrow0.$
Since $\langle \Phi'_\mu(z_n),z_n\rangle=o_n(1)$, we know
$\|z_n\|_{1,\mu}\rightarrow0.$  On the other hand, for any $u,v\in H^1(\mathbb{R}^2)$, in view of (\ref{2.1.0}), the Young inequality and Trudinger-Moser inequality, for any $\delta, \tau>0$, there exist $C_\delta, C_\tau>0$ such that
$$\aligned
\int_{\mathbb{R}^2}|f(u)v|&
\leq\delta|u||v|+C_\delta\int_{\mathbb{R}^2}|u|^{q-1}(e^{\alpha u^2}-1)|v|\\
&\leq\delta|u|^2_2+\delta|v|^2_2+C_\delta\Bigl[\tau|v|^2_2+C_\tau
\int_{\mathbb{R}^2}|u|^{2(q-1)}(e^{2\alpha u^2}-1)\Bigr]\\
&\leq\delta|u|^2_2+\delta|v|^2_2+C_\delta\tau|v|^2_2
+C_\delta C_\tau|u|^{2(q-1)}_{2q}\Bigl(\int_{\mathbb{R}^2}(e^{2q\alpha u^2}-1)\Bigr)^{\frac1q}\\
&\leq\delta|u|^2_2+\delta|v|^2_2+C_\delta\tau|v|^2_2+C_\delta C_\tau|u|^{2(q-1)}_{2q}.
\endaligned$$
Choosing $\delta={\mu}/{4}$ and $\tau={\mu}/(4C_{\frac{\mu}{4}})$ we get $$\int_{\mathbb{R}^2}|f(u)v|\leq\frac{\mu}4|u|^2_2+\frac{\mu}2|v|^2_2+C|u|^{2(q-1)}_{2q}.$$
Similarly for $g$.
Then for any $z=(u,v)\in \mathcal{K}_\mu$, for some $q>2$ we have
 $$0\geq\|u\|^2_\mu+\|v\|^2_\mu
 -C|u|^{2(q-1)}_{2q}-C|v|^{2(q-1)}_{2q}.$$
 Therefore $\|u\|^{2q-4}_\mu+\|v\|^{2q-4}_\mu\geq C$ contradicting $\|z_n\|_{1,\mu}\rightarrow0.$
 Thus $\{z_n\}$ is non-vanishing.
Since $\Phi_\mu$ is $\mathbb{Z}^2$-invariant, up to a translation, we can assume $z_n\rightharpoonup z=(u,v)\neq(0,0)$. It is easy to see that $\Phi'_\mu$ is weakly sequentially continuous, then $\Phi'_\mu(z)=0$ and from Fatou lemma we get $\Phi_\mu(z)\leq \bar{c}_\mu$.

We next show that $c_\mu$ is attained. Assume that $\{\tilde{z}_n\}\subset \mathcal{K}_\mu\backslash\{0\}$ such that $\Phi_\mu(\tilde{z}_n)\rightarrow c_\mu$. Clearly $\{\tilde{z}_n\}$ is a (PS)$_{c_\mu}$ sequence and as above, $\{\tilde{z}_n\}$  is bounded in $E$ and nonvanishing. We may assume that $\tilde{z}_n\rightharpoonup \tilde{z}_0\neq(0,0)$ in $E$.
Taking standard arguments, from Fatou lemma it follows that
$c_\mu= \Phi_\mu(\tilde{z}_0)$.
Then $c_\mu$ is attained and $c_\mu>0$ using (\ref{1.2}).
Noting that $c_\mu$ is also the least energy of $J_\mu$, it is easy to check that $c_\mu\leq \bar{c}_\mu\leq \bar{\bar{{c}}}_\mu$. To prove $\bar{\bar{{c}}}_\mu\leq c_\mu$. As Lemma \ref{l3.7}, we can shows that, if $u\in H^1(\mathbb{R}^2)\backslash\{0\}$ satisfies $\langle J'_\mu(u),u\rangle=0$ then
$J''_\mu(u)[u,u]<0$. Therefore, for any $u\in H^1(\mathbb{R}^2)\backslash\{0\}$, the function $t\mapsto J_\mu(tu)$ has a unique  maximum point $t=t(u)>0$. Denote
$$\mathcal{N}_\mu:=\bigl\{t(u)u: u\in H^1(\mathbb{R}^2)\backslash\{0\}\bigr\},$$
we have $\mathcal{N}_\mu\neq\emptyset$ since $c_\mu$ is attained. Observe that
$\bar{\bar{{c}}}_\mu=\inf_{u\in \mathcal{N}_\mu}J_\mu(u).$
Then $\bar{\bar{{c}}}_\mu\leq c_\mu$ and so $\bar{\bar{{c}}}_\mu=\bar{c}_\mu={c}_\mu$. \ \ \ \ \ $\Box$

To investigate the properties of solutions of (\ref{2.3}), we require in addition the following results of ground states for (\ref{5.1}) whose proof is similar as in  \cite{ZhangCPDE} and we omit.
\begin{lemma}\label{l5.3} Let (H$_1$)-(H$_4$) hold. If in addition, replacing $f$ and $g$ with their odd extensions, then for any ground state $(u,v)$ of the equation (\ref{5.1}), one has $u,v\in C^2(\mathbb{R}^2)$ and $uv>0$ in $\mathbb{R}^2$. Moreover, there exists some point $x_0\in\mathbb{R}^2$
such that $u,v$ are radially
symmetric with respect to the same point $x_0$, namely, $u(x)=u(|x-x_0|)$, $v(x)=v(|x-x_0|)$
and setting $r=|x-x_0|$, without loss of generality, assume that $u,v>0$, one has for $r>0$,
$\frac{\partial u}{\partial r}<0$ and $\frac{\partial v}{\partial r}<0$
as well as
$\Delta u(x_0)<0$, $\Delta v(x_0)<0.$
Moreover, there exist $C,c>0$, such that
$$|D^\alpha u(x)|+|D^\alpha v(x)|\leq C\exp(-c|x-x_0|),\quad x\in\mathbb{R}^2,\ |\alpha|=0,1.$$
\end{lemma}

{\bf 3.1.2. The asymptotically linear case.}
Choose $u\in C^\infty_0(\mathbb{R}^2)$ such that $|u|_2=1$. Set $u_t(x):=tu(tx)$, $t>0$. Then $|u_t|_2=1$ and $|\nabla u_t|_2\rightarrow0$ as $t\rightarrow0$. Hence, by (H$'_3$) there exists $t_0>0$ such that $$|\nabla u_{t_0}|^2_2+\mu|u_{t_0}|^2_2-l_0|u_{t_0}|^2_2\leq|\nabla u_{t_0}|^2_2+|V|_\infty-l_0<0,$$
where $l_0$ is given in (H$'_3$)-(i).
Let \begin{equation}\label{3.2.6}e=\frac{1}{\sqrt{2}}\Bigl(\frac{u_{t_0}}{\|u_{t_0}\|_\mu}, \frac{u_{t_0}}{\|u_{t_0}\|_\mu}\Bigr)\in E^+,\end{equation} and $E_e=E^-\oplus\mathbb{R}e$, we have the following result.

\begin{lemma}\label{l6.2}(1) If $z\in E_e$ and $\|z\|_{1,\mu}\rightarrow+\infty$, then $\Phi_\mu(z)\rightarrow-\infty$  and $\sup_{z\in E_e} \Phi_\mu(z)<+\infty$.\\
\noindent(2) For any $u\in H^1(\mathbb{R}^2)\backslash\{0\}$, then either $J_\mu(tu)\rightarrow+\infty$ or $J_\mu(tu)\rightarrow-\infty$ as $t\rightarrow+\infty$.\end{lemma}
{\bf Proof}: (1) Argue by contradiction we may assume there exists $z_n\in E_e$ satisfying $\|z_n\|_{1,\mu}\rightarrow+\infty$ and $\Phi_\mu(z_n)\geq {-M>-\infty}$. Set
$w_n=\frac{z_n}{\|z_n\|_{1,\mu}}$, $z_n=(z^1_n,z^2_n)$ and $w_n=(w^1_n,w^2_n)$. We assume $w_n\rightharpoonup w=(w^1,w^2)$ in $E$, $w^-_n\rightharpoonup w^-$ in $E$ and $w^+_n\rightarrow w^+\in \mathbb{R}e$. Note that
\begin{equation*}{o_n(1)}\leq \frac{\Phi_\mu(z_n)}{\|z_n\|^2_{1,\mu}}=\frac12(\|w^+_n\|^2_{1,\mu}
-\|w^-_n\|^2_{1,\mu})
-\int_{\mathbb{R}^2}\frac{F(z^1_n)}{\|z_n\|^2_{1,\mu}}-
\int_{\mathbb{R}^2}\frac{G(z^2_n)}{\|z_n\|^2_{1,\mu}}.\end{equation*}
Then $\|w^+_n\|_{1,\mu}\geq\|w^-_n\|_{1,\mu} {+o_n(1)}$. Since $\|w_n\|_{1,\mu}=1$, we have $\|w^+_n\|^2_{1,\mu}\geq\frac12 {+o_n(1)}$. Then $w^+\neq(0,0)$ and so $\frac{w^1+w^2}{2}\neq0$.
Define
$$R_1(t)=F(t)-({l_0}/2)t^2, R_2(t)=G(t)-({l_0}/2)t^2,\quad \forall t\in\mathbb{R}.$$
Since $F(t), G(t)\leq C|t|^2$, we get $R_1(t), R_2(t)\leq C|t|^2$ for any $t\in\mathbb{R}$ and ${R_1(t)}/{t^2}\rightarrow0$, $ {R_2(t)}/{t^2}\rightarrow0$ as $t\rightarrow+\infty.$
Thanks to $w^+\in\mathbb{R}e$, we have $\|w^+\|^2_{1,\mu}<l_0|w^+|^2_2$. Then
$\|w^+\|^2_{1,\mu}< {l_0}|w^1|^2_2+{l_0}|w^2|^2_2.$
Thus there exist bounded domains $\Omega_1,\ \Omega_2\subset\mathbb{R}^2$ such that
\begin{equation}\label{3.2.8}\mathcal{I}:=\frac12(\|w^+\|^2_{1,\mu}-\|w^-\|^2_{1,\mu})-\frac{{l_0}}{2}
\int_{\Omega_1}(w^1)^2
-\frac{{l_0}}{2}\int_{\Omega_2}(w^2)^2<0.\end{equation}
It follows from Lebesgue dominated convergence theorem that
$$\int_{\Omega_1}\frac{R_1(z^1_n)}{\|z_n\|^2_{1,\mu}}=
\int_{\Omega_1}\frac{R_1(z^1_n)|w^1_n|^2}{|z^1_n|^2}\rightarrow0,\quad
\int_{\Omega_2}\frac{R_2(z^2_n)}{\|z_n\|^2_{1,\mu}}=
\int_{\Omega_2}\frac{R_1(z^2_n)|w^2_n|^2}{|z^2_n|^2}\rightarrow0.$$
Then
$$\aligned {o_n(1)}\leq\frac{\Phi_\mu(z_n)}{\|z_n\|^2_{1,\mu}}\leq&\frac{\|w^+_n\|^2_{1,\mu}-\|w^-_n\|^2_{1,\mu}}{2}-
\int_{\Omega_1}\Bigl(\frac{R_1(z^1_n)}{\|z_n\|^2_{1,\mu}}
+\frac{l_0(w^1_n)^2}2\Bigr)\\&
-\int_{\Omega_2}\Bigl(\frac{R_2(z^2_n)}{\|z_n\|^2_{1,\mu}}
+\frac{l_0(w^2_n)^2}2\Bigr)\leq\mathcal{I}+o_n(1),
\endaligned$$
where $\mathcal{I}$ is given in (\ref{3.2.8}). Due to $\mathcal{I}<0$, this is a contradiction. Then $\Phi_\mu(z)\rightarrow-\infty$  and so $\sup_{z\in E_e} \Phi_\mu(z)<+\infty$.

(2) Taking similar arguments of Lemma \ref{l3.6} (1)-(ii),  the conclusion (2) yields. \ \ \ \ \ \ \ \ $\Box$

\begin{lemma}\label{l3.4} Let (H$_1$), (H$_2$) and (H$'_3$) hold. Then

\noindent(1) $\mathcal{K}_\mu\setminus\{0\}\neq\emptyset$ and $c_\mu>0$ is achieved;\\
\noindent(2) $c_\mu=\bar{c}_\mu=\bar{\bar{{c}}}_\mu$, where $\bar{c}_\mu$ and $\bar{\bar{{c}}}_\mu$ are defined as in Lemma \ref{l5.1}.
\end{lemma}
{\bf Proof}: By Lemma \ref{l6.2}, we know $J_\mu(te)\rightarrow-\infty$ as $t\rightarrow\infty$, where $e$ is given in (\ref{3.2.6}). Then the minimax value $\bar{c}_\mu$ is well defined.
 Let $u\in H^1(\mathbb{R}^2)\backslash\{0\}$ we find that the function $t\mapsto J_\mu(tu)$ has at most one nontrivial critical point $t=t(u)>0$. Hence, if we denote
$$\mathcal{N}_\mu:=\bigl\{t(u)u: u\in H^1(\mathbb{R}^2)\backslash\{0\}, t(u)<+\infty\bigr\},$$ then
$\bar{\bar{{c}}}_\mu=\inf_{z\in \mathcal{N}_\mu}J_\mu(z).$ In view of Lemma \ref{l6.1} and Lemma \ref{l2.5}-(3),
as in the proof of  Lemma \ref{l5.1}, we can show that $\mathcal{K}_\mu\setminus\{0\}\neq\emptyset$, $c_\mu>0$ is achieved and $c_\mu=\bar{c}_\mu=\bar{\bar{{c}}}_\mu$.
\ \ \ \ \ \ \ $\Box$

\subsection{Another autonomous system}
For both the superlinear and asymptotically linear cases, we also introduce another autonomous system. If $(u,v)$ is the solution of (\ref{5.1}) with $\mu=V_0$, then
$$(\tilde{u}(x),\tilde{v}(x)):=\Bigl({u}(\sqrt{{V(x_0)}/{V_0}}x),
{v}(\sqrt{{V(x_0)}/{V_0}}x)\Bigr), \ \ \forall x_0\in \mathbb{R}^2,$$ is a solution of
\begin{equation*}\aligned
\left\{ \begin{array}{lll}
-\Delta \tilde{u}+V(x_0)\tilde{u}=\frac{V(x_0)}{V_0}g(\tilde{v})\ & \text{in}\quad \mathbb{R}^2,\\
-\Delta \tilde{v}+V(x_0)\tilde{v}=\frac{V(x_0)}{V_0}f(\tilde{u})\ & \text{in}\quad \mathbb{R}^2.
\end{array}\right.\endaligned
\end{equation*}
Moreover, $\Phi_{V_0}(u,v)=\tilde{\Phi}_{V(x_0)}(\tilde{u},\tilde{v})$, where $\tilde{\Phi}_{V(x_0)}$ is the associated functional of the above system given by
\begin{equation}\label{bu3}\tilde{\Phi}_{V(x_0)}(u,v)=\int_{\mathbb{R}^2}\bigl(\nabla u\nabla v+V(x_0)uv\bigr)-\frac{V(x_0)}{V_0}\int_{\mathbb{R}^2}\bigl(G(v)+F(u)),\ \forall (u,v)\in E.\end{equation}
Then
\begin{equation}\label{3.7.0}\tilde{c}_{V(x_0)}:=
\inf\bigl\{\tilde{\Phi}_{V(x_0)}(z): z\in E\backslash\{(0,0)\}, \ \tilde{\Phi}'_{V(x_0)}(z)=0\bigr\}=c_{V_0}.\end{equation}
As before,  for any $u\in H^1(\mathbb{R}^2)$, there exists $\tilde{h}_{V(x_0)}(u)\in H^1(\mathbb{R}^2)$ such that
\begin{equation}\label{3.7.3} \tilde{\Phi}_{V(x_0)}(u+\tilde{h}_{V(x_0)}(u),u-\tilde{h}_{V(x_0)}(u))
=\max_{v\in H^1(\mathbb{R}^2)}\tilde{\Phi}_{V(x_0)}(u+v,u-v),\end{equation}
and define
\begin{equation}\label{3.7.1}\tilde{J}_{V(x_0)}: H^1(\mathbb{R}^2)\rightarrow\mathbb{R},\ \tilde{J}_{V(x_0)}(u)=\tilde{\Phi}_{V(x_0)}\bigl(u+\tilde{h}_{V(x_0)}(u),
u-\tilde{h}_{V(x_0)}(u)\bigr).\end{equation}

In the case that $f,g$ are superlinear or asymptotically linear at infinity, as in Lemma \ref{l5.1} or Lemma \ref{l3.4}, from (\ref{3.7.0}) we infer
\begin{equation}\label{3.7.2}{c}_{V_0}=\tilde{{c}}_{V(x_0)}=\inf_{u\in H^1(\mathbb{R}^2)\backslash\{0\}}\max_{t>0}\tilde{J}_{V(x_0)}(tu).\end{equation}
Similar to Lemma \ref{l6.2}-(2), we have the following result.
\begin{lemma}\label{l3.10}
For the asymptotically linear case, if $\sup_{t>0}\tilde{J}_{V(x_0)}(tu)\leq M<+\infty$ for some $u\in H^1(\mathbb{R}^2)\backslash\{0\}$, then $\tilde{J}_{V(x_0)}(tu)\rightarrow-\infty$ as $t\rightarrow+\infty$.
\end{lemma}

\section{Proof of Theorems 1.1 and 1.2}
\renewcommand{\theequation}{4.\arabic{equation}}

To show that $J_\epsilon$ has mountain-pass structure, by Lemma \ref{l2.6}, it suffices to show that, there exists $u_0\in H^1(\mathbb{R}^2)$ independent on $\epsilon$, such that $\|u_0\|>r$ and $J_\epsilon(u_0)<0$. For this, we shall investigate the relationship between system (\ref{2.3}) and the limit system (\ref{1.1.1}). In particular,  we need to study the relation between $\bar{h}_\epsilon(u)$ and $\bar{h}_0(u):=\bar{h}_{V_0}(u)$.

\begin{lemma}\label{l3.2} For any $u\in H^1(\mathbb{R}^2)$, $\bar{h}_\epsilon(u)\rightarrow\bar{h}_{0}(u)$  in $H^1(\mathbb{R}^2)$ as $\epsilon\rightarrow0$.\end{lemma}
{\bf Proof}: For any sequence $\epsilon_n\rightarrow 0^+$, put $v_n:=\bar{h}_{\epsilon_n}(u)$ and $v_0:=\bar{h}_{0}(u)$. We shall prove that
$v_n\rightarrow v_0$ in $H^1(\R^2)$.
Firstly, using (\ref{3.1.3}), it is easy to see that $\{v_n\}\subset H^1(\R^2)$ is bounded. Up to a subsequence, we may assume
$v_n\rightharpoonup v^*$ {in} $H^1(\R^2)$
and thus $u\pm v_n\rightharpoonup u\pm v^*$ in $H^1(\R^2)$.
Recalling that
\begin{equation*}\langle\Phi'_{\epsilon_n}(u+v_n,
u-v_n),(\varphi,-\varphi)\rangle=0, \quad \forall \varphi\in H^1(\mathbb{R}^2),\end{equation*}
that is, for any $\varphi\in H^1(\mathbb{R}^2)$,
\begin{equation}\label{eq:20220410-e5}
\int_{\R^2}[-2\nabla v_n \nabla \varphi-2V(\epsilon_n x)v_n\varphi+\bar{g}(\epsilon_n x,u-v_n)\varphi -\bar{f}(\epsilon_n x, u+v_n)\varphi]=0.\end{equation}
It follows that
\begin{equation*}
\int_{\R^2}[-2\nabla v^* \nabla \varphi-2V_0v^*\varphi+g(u-v_0)\varphi -f( u+v_0)\varphi]=0, \forall \varphi\in H^1(\R^2).
\end{equation*}
By the definition of $\bar{h}_{0}$ (see (\ref{3.3})), we obtain $v^*=\bar{h}_{0}(u)=v_0$.
Testing \eqref{eq:20220410-e5} by $\varphi=v_n$, there holds
 \begin{equation*}\label{eq:20220410-e7}
 2\|v_n\|_{\epsilon_n}^{2}=\int_{\R^2}[\bar{g}(\epsilon_n x,u-v_n)v_n-\bar{f}(\epsilon_n x, u+v_n) v_n ].
 \end{equation*}
 Similarly, we have
 \begin{equation*}\label{eq:20220410-e8}
 2\|v_0\|_{0}^{2}=\int_{\R^2}[g(u-v_0)v_0-f(u+v_0) v_0 ].
 \end{equation*}
Therefore
 \begin{equation}\label{eq:20220411-e0}
 \aligned
 2\|v_n-v_0\|_{\epsilon_n}^{2}=&2\left(\|v_n\|_{\epsilon_n}^{2}
 -\|v_0\|_{\epsilon_n}^{2}\right)+o_n(1)
 =2\left(\|v_n\|_{\epsilon_n}^{2}-\|v_0\|_{0}^{2}\right)+o_n(1)\\
 =&o_n(1)+\int_{\R^2}\left[\bar{g}(\epsilon_n x, u-v_n)v_n -g(u-v_0)v_0\right]\\
 &-\int_{\R^2}\left[\bar{f}(\epsilon_n x, u+v_n)v_n -f(u+v_0)v_0\right].
 \endaligned
 \end{equation}
 Noting that
 \begin{equation}\label{eq:20220411-e1}
 \aligned
 &\int_{\R^2}\left[\bar{f}(\epsilon_n x, u+v_n)v_n -f(u+v_0)v_0\right]\\
 =&\int_{\R^2}\left[\bar{f}(\epsilon_n x, u+v_n)-\bar{f}(\epsilon_n x, u+v_0)\right](v_n-v_0)+\int_{\R^2}\bar{f}(\epsilon_n x, u+v_0)(v_n-v_0)\\
 &+\int_{\R^2}\bar{f}(\epsilon_n x, u+v_n)v_0 -\int_{\R^2}f(u+v_0)v_0.
 \endaligned
 \end{equation}
Using $(H''_4)$ one easily has
\begin{equation}
\int_{\R^2}\left[\bar{f}(\epsilon_n x, u+v_n)-\bar{f}(\epsilon_n x, u+v_0)\right](v_n-v_0) \geq 0.
\end{equation}
Since $v_n\rightharpoonup v_0$ in $H^1(\mathbb{R}^2)$, we infer
\begin{equation}\label{eq:20220411-e2}
\aligned
&\int_{\R^2}\bar{f}(\epsilon_n x, u+v_0)(v_n-v_0)\\
=&\int_{\Lambda_{\epsilon_n}}f(u+v_0)(v_n-v_0) +\int_{\Lambda_{\epsilon_n}^{c}}\bar{f}(\epsilon_n x, u+v_0)(v_n-v_0)\\
=&o_n(1)+\int_{\Lambda_{\epsilon_n}^{c}}\bar{f}(\epsilon_n x, u+v_0)(v_n-v_0)
\endaligned
\end{equation}
 and
 \begin{equation}\label{eq:20220411-e3}
 \aligned
 &\int_{\R^2}\bar{f}(\epsilon_n x, u+v_n)v_0 -\int_{\R^2}f(u+v_0)v_0\\
 =&\int_{\Lambda_{\epsilon_n}} [f(u+v_n)-f(u+v_0)]v_0
 +\int_{\Lambda_{\epsilon_n}^{c}} [\bar{f}(\epsilon_n x, u+v_n) -f(u+v_0)]v_0\\
 =&o_n(1).
 \endaligned
 \end{equation}
In view of \eqref{eq:20220411-e1}-\eqref{eq:20220411-e3}, we obtain
 \begin{equation}\label{eq:20220411-e4}
 \aligned
 &\int_{\R^2}\left[\bar{f}(\epsilon_n x, u+v_n)v_n -f(u+v_0)v_0\right]
\geq\int_{\Lambda_{\epsilon_n}^{c}}\bar{f}(\epsilon_n x, u+v_0)(v_n-v_0)+o_n(1).
 \endaligned
 \end{equation}
 Similarly, there holds
 \begin{equation}\label{eq:20220411-e5}
 \aligned
 &\int_{\R^2}\left[\bar{g}(\epsilon_n x, u-v_n)v_n -g(u-v_0)v_0\right]
\leq \int_{\Lambda_{\epsilon_n}^{c}}\bar{g}(\epsilon_n x, u-v_0)(v_n-v_0)+o_n(1).
 \endaligned
 \end{equation}
From $(H''_3)$, \eqref{eq:20220411-e0},\eqref{eq:20220411-e4} and \eqref{eq:20220411-e5}, it follows that
 \begin{equation*}
 \aligned
 2\|v_n-v_0\|_{\epsilon_n}^{2}\leq&o_n(1)+\int_{\Lambda_{\epsilon_n}^{c}}
 \bigl(\bar{g}(\epsilon_n x, u-v_0)-\bar{f}(\epsilon_n x, u+v_0)\bigr)(v_n-v_0)\\
 \leq& o_n(1)+ \frac{\inf_{\R^2}V}{2}  \bigl(\|u-v_0\|_{L^2(\Lambda_{\epsilon_n}^{c})}+\|u+v_0\|_{L^2(\Lambda_{\epsilon_n}^{c})} \bigr)|v_n-v_0|_{2}\\
 \leq&o_n(1)+o_n(1)\|v_n-v_0\|_{\epsilon_n},
 \endaligned
 \end{equation*}
 which implies $v_n\rightarrow v_0$ in $H^1(\R^2)$.
\ \ \  \ $\Box$

\begin{lemma}\label{l3.3}For $\epsilon>0$ small enough, there exists $u_0\in H^1(\mathbb{R}^2)$ (independent of $\epsilon$) such that $\|u_0\|_\epsilon>r$, where $r$ is given in Lemma \ref{l2.6}, and $J_\epsilon(u_0)<0$.\end{lemma}
{\bf Proof}: Let $(v_1,v_2)$ be a ground state of the limit problem (\ref{5.1}) with $\mu=V_0$. Then $J_{V_0}(u)=c_{V_0}$ where $u=\frac{v_1+v_2}{2}$ and
$$c_{V_0}=\inf_{u\in H^1(\mathbb{R}^2)\backslash\{0\}}\max_{t>0}J_{V_0}(tu).$$
For the asymptotically linear case, observe that  $\sup_{t>0}J_{V_0}(tu)=J_{V_0}(u)=c_{V_0}$, then Lemma \ref{l6.2}-(2) implies that $J_{V_0}(tu)\rightarrow-\infty$ as $t\rightarrow+\infty$. Combining with Lemma \ref{l4.7}-(2) we know there exists $t_0>0$ large enough such that
\begin{equation*}\aligned J_{V_0}(t_0u)=&t^2_0\|u\|^2_0-\|\bar{h}_0(t_0u)\|^2_0-
\int_{\mathbb{R}^2}[F(t_0u+\bar{h}_0(t_0u))+
G(t_0u-\bar{h}_0(t_0u))]<-1.
 \endaligned\end{equation*}
Hence, there is $R_0>0$ such that
\begin{equation}\label{3.13}\aligned \mathcal{I}:=t^2_0\|u\|^2_0-\|\bar{h}_0(t_0u)\|^2_0
-\int_{B_{R_0}}[F(t_0u+\bar{h}_0(t_0u))+
 G(t_0u-\bar{h}_0(t_0u))]<-\frac12.\endaligned\end{equation}
Recalling that $V(\epsilon x)\rightarrow V_0$ uniformly on bounded sets of $\mathbb{R}^2$, it follows from
Lemma \ref{l3.2} that
$$\aligned J_\epsilon(t_0 u)\leq t^2_0\|u\|^2_\epsilon-\|\bar{h}_\epsilon(t_0u)\|^2_\epsilon
-\int_{\Lambda_\epsilon}
[{F}(t_0u+\bar{h}_\epsilon(t_0u))+{G}(t_0u-
\bar{h}_\epsilon(t_0u))]\leq\mathcal{I}+o_\epsilon(1),
\endaligned$$
 as $\epsilon\rightarrow0$.
 Therefore, by (\ref{3.13}) there is $\epsilon_0>0$ such that $J_\epsilon(t_0u)<0$ for all $0<\epsilon<\epsilon_0$.  \ \ \ \ \ $\Box$

\begin{lemma}\label{l3.5.0} For $\epsilon$ small enough, the system (\ref{2.3}) has a nontrivial solution.\end{lemma}
{\bf Proof}: For the superlinear case, by Lemma \ref{l2.2} and Lemma \ref{l2.5}-(2),  it is easy to see that $J_\epsilon$ satisfies the (PS) condition. In view of Lemmas \ref{l2.6} and \ref{l3.3}, applying the mountain-pass theorem with (PS) condition, we know
\begin{equation}\label{3.18.0}c_\epsilon=\inf_{v\in \Gamma_\epsilon}\max_{t\in[0,1]}J_\epsilon(\nu(t)),\end{equation}
where $\Gamma_\epsilon=\{\nu\in C([0,1], H^1(\mathbb{R}^2)):\nu(0)=0, J_\epsilon(\nu(1))<0\}$, is a critical value of $J_\epsilon$, also for $\Phi_\epsilon$ and $\tau\leq c_\epsilon<+\infty$ with $\tau$ given in Lemma \ref{l4.7}-(1).

For the asymptotically linear case,
  from Lemma \ref{l6.1} and Lemma \ref{l2.5}-(3) it follows that $J_\epsilon$ satisfies (Ce) condition.
  In view of Lemmas \ref{l2.6} and \ref{l3.3}, applying the mountain-pass theorem with (Ce) condition, we know $c_\epsilon$ given in (\ref{3.18.0}) is a critical value of $\Phi_\epsilon$.
 \ \ \ \ \ $\Box$

\begin{lemma}\label{l3.8}
(1) For $\epsilon$ small enough, $c_\epsilon=\inf_{u\in H^1(\mathbb{R}^2)\backslash\{0\}}\max_{t\geq0}J_\epsilon(tu)$, where $c_\epsilon$ is given in (\ref{3.18.0}).\\
\noindent(2) $c_\epsilon\leq c_{V_0}+o_\epsilon(1)$ as $\epsilon\rightarrow0$.
\end{lemma}
{\bf Proof}: (1) Setting $d_\epsilon=\inf_{u\in H^1(\mathbb{R}^2)\backslash\{0\}}\max_{t\geq0}J_\epsilon(tu)$. By Lemma \ref{l3.6}, we get $d_\epsilon\geq c_\epsilon$. It suffices to show the other inequality. Letting $u\in H^1(\mathbb{R}^2)$, from Lemma \ref{l3.7} it follows that the function $t\mapsto J_\epsilon(tu)$ has at most one nontrivial critical point $t=t(u)>0$. Denote
$$\mathcal{N}_\epsilon:=\{t(u)u: u\in H^1(\mathbb{R}^2)\backslash\{0\},t(u)<+\infty\}.$$
By Lemma \ref{l3.3}, we know $\mathcal{N}_\epsilon\neq\emptyset$. Moreover, it is easy to verify that
$d_\epsilon=\inf_{u\in \mathcal{N}_\epsilon}J_\epsilon(u).$
So we only need to show that given $\nu\in \Gamma_\epsilon$, there exists $t_0\in [0,1]$ such that $\nu(t_0)\in \mathcal{N}_\epsilon$. Otherwise, we assume that $\nu([0,1])\cap \mathcal{N}_\epsilon=\emptyset$. Similar to the proof of Lemma \ref{l2.6}, we have
$\langle J'_\epsilon(\nu(t)),\nu(t)\rangle>0$ for small $t>0$. Since the function $t\rightarrow \langle J'_\epsilon(\nu(t)),\nu(t)\rangle$ is continuous and $\langle J'_\epsilon(\nu(t)),\nu(t)\rangle\neq0$ for all $t\in(0,1]$, we get
$\langle J'_\epsilon(\nu(t)),\nu(t)\rangle>0$ for all $t\in(0,1].$
Then
$ J_\epsilon(\nu(t)){\ge\frac12 \langle J'_\epsilon(\nu(t)),\nu(t) \rangle >0}$, for all $t\in[0,1]$, which contradicts the definition of $\Gamma_\epsilon$.  So $\nu(t)$ crosses
$\mathcal{N}_\epsilon$ provided $\nu\in \Gamma_\epsilon$.

(2) Let $(v_1,v_2)$ be a ground state of the limit problem (\ref{5.1}) with $\mu=V_0$. Then $w\in \mathcal{N}_{V_0}$ and $J_{V_0}(w)=c_{V_0}$ with $w=\frac{v_1+v_2}{2}$.  Similar to the argument of Lemma \ref{l3.3}, there exists $t_0>0$ large enough such that $J_\epsilon(t_0 w)\leq-\frac12$ with small enough $\epsilon$. Then there exists $t_{w}\in (0,1)$ such that $t_{w}t_0w\in \mathcal{N}_\epsilon$. Note that $c_\epsilon=\inf_{\mathcal{N}_\epsilon}J_\epsilon$ in the conclusion (1), we have
$c_\epsilon\leq J_\epsilon(t_{w}t_0w)$. Define the family $\{\Upsilon_\epsilon\}\subset C([0,t_0])$ by
$\Upsilon_\epsilon(t)=J_\epsilon(tw)-J_{V_0}(tw).$
By the boundedness of $\bar{h}_\epsilon$ and $\bar{h}'_\epsilon$, one easily conclude that $\Upsilon_\epsilon$ and $\Upsilon'_\epsilon$ are uniformly bounded. Then $\Upsilon'_\epsilon$ is equicontinuous. Applying the  Arzel$\grave{a}$-Ascoli theorem to $\Upsilon_\epsilon$, we deduce
\begin{equation*}J_\epsilon(tw)=J_{V_0}(tw)+o_\epsilon(1)\ \ \text{uniformly in } \ t\in[0,t_0],\end{equation*}
as $\epsilon\rightarrow0$. Consequently, $c_\epsilon\leq c_{V_0}+o_\epsilon(1)$ as $\epsilon\rightarrow0$. This ends the proof. \ \ \ \ $\Box$

To show the family of solutions converge strongly as $\epsilon\rightarrow0$, we need the following lemma.
\begin{lemma}\label{l4.1}
Assume $u_n\rightharpoonup u$ in $H^1(\mathbb{R}^2)$. Then
\begin{equation}\label{4.5}\int_{\mathbb{R}^2}F(u_n)
-\int_{\mathbb{R}^2}F(u_n-u)-\int_{\mathbb{R}^2}F(u)=o_n(1),\end{equation}
\begin{equation}\label{4.6}\int_{\mathbb{R}^2}f(u_n)\varphi
-\int_{\mathbb{R}^2}f(u_n-u)\varphi-\int_{\mathbb{R}^2}f(u)\varphi
=o_n(1)\|\varphi\|_\epsilon,\end{equation}
uniformly in $\varphi\in H^1(\mathbb{R}^2)$ with $\|\varphi\|_\epsilon\leq1$, and similar results hold for  $G$ and $g$.
\end{lemma}
{\bf Proof}: Firstly we show (\ref{4.5}). Set $v_n=u_n-u$, we know
$F(u_n)-F(v_n)=f(v_n+tu)u$ with $0<t:=t(n,x)<1$. Set $w_n=v_n+tu$, and assume $\|w_n\|^2_\epsilon\leq M$. By (\ref{2.1.0}), for some $\alpha<\frac{2\pi}{M}$, for any $\delta>0$, there exists $C_\delta>0$ such that
$$\aligned |F(v_n+u)-F(v_n)|&=
|f(v_n+tu)u|\leq  \delta|v_n||u|+\delta|u|^2+
C_\delta|u|H_n,\endaligned$$
where $H_n=e^{\alpha w^2_n}-1$. Using the Young inequality, we have
\begin{equation}\label{4.7}\aligned |F(v_n+u)-F(v_n)-F(u)|\leq & \delta|v_n|^2+C_1|u|^2+\delta H^{2}_n+C|u|H,\endaligned\end{equation}
where $H=e^{\alpha u^2}-1$. By the Trudinger-Moser inequality, we get
\begin{equation*}\int_{\mathbb{R}^2}H^{2}_n
\leq \int_{\mathbb{R}^2}[e^{2\alpha\|w_n\|^2_\epsilon\frac{w^2_n}{\|w_n\|^2_\epsilon}}-1]\leq C.\end{equation*}
Define
$$\psi_{\delta,n}=\max\bigl\{|F(v_n+u)-F(v_n)-F(u)|-\delta|v_n|^2
-\delta H^{2}_n,0\bigr\}.$$
From (\ref{4.7}) it is easy to see
$\psi_{\delta,n}\rightarrow 0$ a.e. in $\mathbb{R}^2$, and $\psi_{\delta,n}\leq C_1|u|^2+C|u|H\in L^1(\mathbb{R}^2).$
The Lebesgue theorem implies that
$\int_{\mathbb{R}^2}\psi_{\delta,n}\rightarrow0,\ \ \text{as}\ n\rightarrow\infty.$
Thus
$$\limsup_{n\rightarrow\infty}\int_{\mathbb{R}^2}|F(v_n+u)-F(v_n)-F(u)|
\leq\limsup_{n\rightarrow\infty}\int_{\mathbb{R}^2}
(\delta|v_n|^2
+\delta H^{2}_n+\psi_{\delta,n})
\leq C\delta.$$
By the arbitrariness of $\delta$,
(\ref{4.5}) yields. Using (\ref{2.2.0}), (\ref{4.6}) can be deduced similarly. \ \ \ \ $\Box$

\begin{lemma} \label{l4.2} Let $z_{n}=(u_n,v_n)$ be nontrivial solutions of (\ref{2.3}) obtained in Lemma \ref{l3.5.0} with
$\epsilon_n\rightarrow0$. If in addition $V$ is uniformly continuous, then there is $x_{n}\in\mathbb{R}^2$
such that $dist(\epsilon_n x_{n},\mathcal{V})\rightarrow0$ with $\mathcal{V}=\{x\in \Lambda: V(x)=V_0\}$ and the
sequence $\bar{z}_{n}(x):=z_{n}(x+x_n)$ converges strongly in $E$ to a ground
state $z_0=(u_0,v_0)$ of (\ref{1.1.1}).
 \end{lemma}
{\bf Proof}: By Lemma \ref{l3.8}, we know $c_{\epsilon_n}\leq c_{V_0}+o_n(1)$, then as the proof of Lemmas \ref{l2.2} and \ref{l6.1}, we know $\{z_n\}$ is bounded in $E$ and assume $\|z_n\|^2_{1,\epsilon}\leq M$ for some constant $M>0$. If $\{z_n\}$ is vanishing, i.e.
 $$\lim_{n\rightarrow\infty}\sup_{y\in\mathbb{R}^2}\int_{B_R(y)}
 (u^2_n+v^2_n)=0, \ \ \forall R>0,$$
then P.L. Lions compactness lemma implies that $u_n\rightarrow0$ and $v_n\rightarrow0$ in $L^r(\mathbb{R}^2)$ for any $r>2$.
As in the proof of Lemma \ref{l5.1} we have
$\|z_n\|^2_{1,\epsilon}\rightarrow0$. In addition, similar to (\ref{5.3}) we get $\bar{F}(\epsilon x,u_n)\rightarrow0$ and $\bar{G}(\epsilon x,v_n)\rightarrow0$. Then $c_{\epsilon_n}\rightarrow0$, contradicts the fact that $c_{\epsilon_n}\geq \tau>0$.
So $\{z_n\}$ is nonvanishing and there exist $x_n\in\mathbb{R}^2$ and $\delta>0$ such that
\begin{equation}\label{4.9.0}\int_{B_1(x_n)}(u^2_n+v^2_n)\geq\delta.\end{equation}
{\bf Claim 1}: $\epsilon_nx_n\rightarrow x_0\in \Lambda$ and $V(x_0)=V_0$.

Letting $\bar{z}_n=z_n(x+x_n)$, then $\bar{z}_n \rightharpoonup z_0=(u_0,v_0)\neq(0,0)$. Moreover, $\bar{z}_n$ satisfies the system
\begin{equation}\label{4.2}\aligned
\left\{ \begin{array}{lll}
-\Delta u+V(\epsilon_n x+\epsilon_n x_n)u=\bar{g}(\epsilon_n x+\epsilon_nx_n,v)\ & \text{in}\quad \mathbb{R}^2,\\
-\Delta v+V(\epsilon_n x+\epsilon_n x_n)v=\bar{f}(\epsilon_n x+\epsilon_n x_n,u)\ & \text{in}\quad \mathbb{R}^2.
\end{array}\right.\endaligned
\end{equation}
According to (\ref{4.9.0}) we discuss into two cases.

{\bf Case 1:} $\int_{B_1(x_n)}\Bigl[(1-\chi(\epsilon_n x))u^2_n+(1-\chi(\epsilon_n x))v^2_n\Bigr]\geq\frac{\delta}{2}$.
In this case, $B_1(x_n)\cap \{\mathbb{R}^2\setminus supp\chi(\epsilon_n x) \} \neq\emptyset$.
Then $|\epsilon_n x_n|\rightarrow\infty $ or $\epsilon_n x_n\rightarrow x_0\in \mathbb{R}^2\backslash\Lambda$, and further we assume that $x_0\not\in \partial\Lambda$ since we can argue as the following Case 2 if $x_0\in \partial\Lambda$). Observe that $V\in L^\infty(\mathbb{R}^2)$, we may assume $V_1=\lim_{n\rightarrow\infty}V(\epsilon_n x_n)$. Since $V$ is uniformly continuous, we deduce that
$z_0=(u_0,v_0)$ satisfies
\begin{equation*}\aligned
\left\{ \begin{array}{lll}
-\Delta u_0+V_1u_0=\tilde{g}(v_0)\ & \text{in}\quad \mathbb{R}^2,\\
-\Delta v_0+V_1v_0=\tilde{f}(u_0)\ & \text{in}\quad \mathbb{R}^2.
\end{array}\right.\endaligned
\end{equation*}
Then
$$\aligned &\frac12\bigl(\|u_0\|^2_{V_1}+\| v_0|^2_{V_1}\bigr)
=\int_{\mathbb{R}^2}\tilde{g}(v_0)u_0
+\int_{\mathbb{R}^2}\tilde{f}(u_0)v_0
\leq\frac{\inf_{\mathbb{R}^2}V}{2}(|u_0|^2_2+
|v_0|^2_2).\endaligned$$
So $u_0=v_0=0$, contradicts $(u_0,v_0)\neq(0,0)$.

{\bf Case 2}: $\int_{B_1(x_n)}\bigl[\chi(\epsilon_n x)u^2_n+\chi(\epsilon_n x)v^2_n\bigr]\geq\frac{\delta}{2}$. If this case occurs, we have
$B_1(x_n)\cap supp\chi(\epsilon_n x) \neq\emptyset$. Then $\epsilon_n x_n\rightarrow x_0\in \bar{\Lambda}$ as $n\rightarrow\infty$ and $z_0=(u_0,v_0)$ satisfies
\begin{equation}\label{4.3}\aligned
\left\{ \begin{array}{lll}
-\Delta u_0+V(x_0)u_0={g}_\infty(x,v_0)\ & \text{in}\quad \mathbb{R}^2,\\
-\Delta v_0+V(x_0)v_0={f}_\infty(x,u_0)\ & \text{in}\quad \mathbb{R}^2,
\end{array}\right.\endaligned
\end{equation}
where ${g}_\infty(x,s)=\chi_\infty g(s)+(1-\chi_\infty)\tilde{g}(s)$, ${f}_\infty(x,s)=\chi_\infty f(s)+(1-\chi_\infty)\tilde{f}(s)$ and $\chi_\infty$ is either a characteristic
function of a half-space of $\mathbb{R}^2$ provided $\limsup_{\epsilon\rightarrow0}dist(x_n,\partial \Lambda_{\epsilon_n})<+\infty$
or $\chi_\infty\equiv1$, this can be seen by
the fact that $\chi_\Lambda(\epsilon_n(\cdot+x_{n}))\rightarrow\chi_\infty$ a.e. in $\mathbb{R}^2$.
 Denote
the functional of (\ref{4.3}) as
$$\Phi_\infty(u,v)=\int_{\mathbb{R}^2}(\nabla u\nabla v+V(x_0)uv)-\int_{\mathbb{R}^2}[G_\infty(x,v)+F_\infty(x,u)],\quad\forall (u,v)\in E,$$
where $F_{\infty}(x,s)=\int^s_0f_\infty(x,t)dt$, $G_{\infty}(x,s)=\int^s_0g_\infty(x,t)dt$.
By noting that
$$\int_{\mathbb{R}^2}F_\infty(x,u)\leq\int_{\mathbb{R}^2}F(u),\quad \int_{\mathbb{R}^2}G_\infty(x,v)\leq\int_{\mathbb{R}^2}G(v),$$ we have
$$\Phi_\infty(u,v)=\tilde{\Phi}_{V(x_0)}(u,v)+\frac{V(x_0)}{V_0}
\int_{\mathbb{R}^2}(G(v)+F(u))-\int_{\mathbb{R}^2}[F_\infty(x,u)+G_\infty(x,v)],$$
{for all} $(u,v)\in E$, where $\tilde{\Phi}_{V(x_0)}$ is given in (\ref{bu3}). As in Section 2.2, define $h_\infty: H^1(\mathbb{R}^2)\rightarrow H^1(\mathbb{R}^2)$ and $J_\infty: H^1(\mathbb{R}^2)\rightarrow\mathbb{R}$ by
$$J_\infty(u)=\Phi_\infty(u+h_\infty(u),u-h_\infty(u))=\max_{v\in H^1(\mathbb{R}^2)}\Phi_\infty(u+v,u-v).$$
Furthermore, as Lemma \ref{l3.7} we can show that, if $u\in H^1(\mathbb{R}^2)\backslash\{0\}$ satisfies $\langle J'_\infty(u),u\rangle=0$, then $J''_\infty(u)[u,u]<0$. Since we already have that $z_0=(u_0,v_0)\neq(0,0)$ is a critical point of $\Phi_\infty$, we then infer that $\frac{u_0+v_0}{2}$ is a critical point of $J_\infty$ and $J_\infty(\frac{u_0+v_0}{2})=\max_{t>0}J_\infty(t\frac{u_0+v_0}{2})$.
We claim that
\begin{equation}\label{4.2.5} \text{there exists}\ t_0>0\ \text{such that}\ \tilde{J}_{V(x_0)}(t_0 \frac{u_0+v_0}{2})=\max_{t>0}\tilde{J}_{V(x_0)}(t \frac{u_0+v_0}{2}),\end{equation}
where $ \tilde{J}_{V(x_0)}$ is given in (\ref{3.7.1}).
In fact, for the superlinear case, (\ref{4.2.5}) is easily obtained. For the asymptotically linear case, if
\begin{equation}\label{5.15.0}\sup_{t>0}\tilde{J}_{V(x_0)}(t \frac{u_0+v_0}{2})\leq M,  \ \text{for some}\ M>0,\end{equation}
holds,  by Lemma \ref{l3.10} we then obtain that $ \tilde{J}_{V(x_0)}(tu)\rightarrow-\infty$ as $t\rightarrow\infty$ and so (\ref{4.2.5}) yields. So it suffices to show (\ref{5.15.0}).
Actually, for any $t>0$, denote
$$\xi_1(t)=t\frac{u_0+v_0}{2}+\tilde{h}_{V(x_0)}(t\frac{u_0+v_0}{2}),\quad
\xi_2(t)=t\frac{u_0+v_0}{2}-\tilde{h}_{V(x_0)}(t\frac{u_0+v_0}{2}),$$
where $\tilde{h}_{V(x_0)}$ is given in (\ref{3.7.3}).
Then for any $t>0$
\begin{equation}\label{4.2.6}\aligned&\tilde{J}_{V(x_0)}(t\frac{u_0+v_0}{2})=\tilde{\Phi}_{V(x_0)}
(\xi_1(t),\xi_2(t)) \leq  \tilde{\Phi}_{V(x_0)}\bigl(\xi_1(t),
\xi_2(t)\bigr)
\\&+\int_{\mathbb{R}^2}\frac{V(x_0)}{V_0}\bigl(G(\xi_2(t))+F(\xi_1(t))\bigr)
-\int_{\mathbb{R}^2}G_{\infty}(x,\xi_2(t))
-\int_{\mathbb{R}^2}F_{\infty}(x,\xi_1(t))\\
=&\Phi_\infty(\xi_1(t),\xi_2(t))
\leq\Phi_\infty\Bigl(t\frac{u_0+v_0}{2}+{h}_{\infty}(t\frac{u_0+v_0}{2}),
t\frac{u_0+v_0}{2}-{h}_{\infty}(t\frac{u_0+v_0}{2})\Bigr)\\=& J_\infty(t\frac{u_0+v_0}{2})\leq J_\infty(\frac{u_0+v_0}{2})=\Phi_\infty({u_0},{v_0}).
\endaligned\end{equation}
Thus (\ref{4.2.5}) and (\ref{5.15.0}) yield. By  (\ref{3.7.2}), (\ref{4.2.5}) and (\ref{4.2.6}) we deduce
\begin{equation}\label{4.2.8}
c_{V_0}\leq\tilde{J}_{V(x_0)}(t_0\frac{u_0+v_0}{2})
\leq\Phi_\infty({u_0},{v_0}).\end{equation}

On the other hand, from Fatou Lemma it follows that
$$\aligned c_{\epsilon_n}&={\Phi}_{\epsilon_n}(z_n)-\frac12\langle {\Phi}'_{\epsilon_n}(z_n),z_n\rangle=\int_{\mathbb{R}^2}(\hat{F}(\epsilon_n x,u_n)+\hat{G}(\epsilon_n x,v_n))\\&=\int_{\mathbb{R}^2}\bigl[\hat{F}(\epsilon_n x+\epsilon_n x_n,u_n(x+x_n))+\hat{G}(\epsilon_n x+\epsilon x_n,v_n(x+x_n))\bigr]\\
&\geq \int_{\mathbb{R}^2}\hat{F}_\infty( x,u_0)+\int_{\mathbb{R}^2}\hat{G}_\infty( x,v_0)=\Phi_\infty(z_0)-\frac12\langle \Phi'_\infty(z_0),z_0\rangle=\Phi_\infty(u_0,v_0),\endaligned$$
where $\hat{F}_\infty( x,s)=\frac12f_\infty(x,s)s-F_\infty(x,s)$ and
$\hat{G}_\infty( x,s)=\frac12g_\infty(x,s)s-G_\infty(x,s)$ for all $(x,s)\in\mathbb{R}^2\times \mathbb{R}$. Therefore, together with (\ref{4.2.8}), (\ref{4.2.6}) and $c_{\epsilon_n}\leq c_{V_0}+o_n(1)$, we know that $V(x_0)=V_0$, $x_0\in \Lambda$ and $\chi_\infty\equiv1$.

{\bf Claim 2}: $\bar{z}_n\rightarrow z_0$ in $H^1(\mathbb{R}^2)$.

Recall that $\bar{z}_n(x)=(u_n(x+x_n),v_n(x+x_n))$ satisfies the system (\ref{4.2}), whose corresponding functional  is denoted as
$$\aligned\tilde{\Phi}_{\epsilon_n}(u,v)
=&\int_{\mathbb{R}^2}(\nabla u\nabla v+{\tilde{V}}_{\epsilon_n}(x)uv)\\&
-\int_{\mathbb{R}^2}\bar{F}(\epsilon_nx+\epsilon_nx_n, u)-\int_{\mathbb{R}^2}\bar{G}(\epsilon_nx+\epsilon_nx_n, v),\ \ \text{for all}\ (u,v)\in E,\endaligned$$
where $\tilde{V}_{\epsilon_n}(x)=V(\epsilon_n x+\epsilon_nx_n)$.
Moreover, $z_0=(u_0,v_0)$ satisfies the system (\ref{1.1.1})
and
\begin{equation}\label{4.11} \Phi_{0}(z_0)=c_{V_0}
=c_{\epsilon_n}+o_n(1)={\Phi}_{\epsilon_n}(z_{n})+o_n(1)
=\tilde{{\Phi}}_{\epsilon_n}(\bar{z}_{n})+o_n(1),
\end{equation}
where $\Phi_{0}=\Phi_{V_0}$ given before. We argue by contradiction. Otherwise, set $w_{n,1}(x)=z_n(x)-{z_0}(x-x_n)$,
then $w_n(x)=w_{n,1}(x+x_n)\rightharpoonup0$ and $w_{n}\not\rightarrow0$ in $E$. For any $\varphi=(\varphi_1,\varphi_2)\in E$, similar to Lemma \ref{l4.1} we have
$$\tilde{{\Phi}}_{\epsilon_n}(w_n)=\tilde{{\Phi}}_{\epsilon_n}(z_n)
-\tilde{{\Phi}}_{\epsilon_n}(z_0)+o_n(1).$$
$$\langle\tilde{{\Phi}}'_{\epsilon_n}(w_n),\varphi\rangle=\langle\tilde{{\Phi}}'_{\epsilon_n}(z_n),\varphi\rangle
-\langle\tilde{{\Phi}}'_{\epsilon_n}(z_0),\varphi\rangle+o_n(1)\|\varphi\|_E.$$
Using the fact that  $\tilde{V}_{\epsilon_n}(x)\rightarrow V_0$, $\bar{F}(\epsilon_nx+\epsilon_nx_n,s)\rightarrow \bar{F}(x_0,s)=F(s)$ as $n\rightarrow\infty$ uniformly
on any bounded set of $x$, one easily has
\begin{equation*}\aligned &\int_{\mathbb{R}^2}\tilde{V}_{\epsilon_n}(x) u_0\varphi_1(x+x_n)-\int_{\mathbb{R}^2}V_0 u_0\varphi_1(x+x_n)=o_n(1)\|\varphi_1\|_{\epsilon_n},\\
&\int_{\mathbb{R}^2}\tilde{V}_{\epsilon_n}(x) v_0\varphi_2(x+x_n)-\int_{\mathbb{R}^2}V_0 v_0\varphi_2(x+x_n)=o_n(1)\|\varphi_2\|_{\epsilon_n},\endaligned\end{equation*}
and
\begin{equation*}\aligned &\int_{\mathbb{R}^2}\bar{f}({\epsilon_n}x+\epsilon_nx_n,u_0)\varphi_2(x+x_n)
-\int_{\mathbb{R}^2}f(u_0)\varphi_2(x+x_n)=o_n(1)\|\varphi_2\|_{\epsilon_n},\\
&\int_{\mathbb{R}^2}\bar{g}({\epsilon_n}x+\epsilon_nx_n,v_0)\varphi_1(x+x_n)
-\int_{\mathbb{R}^2}g(v_0)\varphi_1(x+x_n)=o_n(1)\|\varphi_1\|_{\epsilon_n}.
\endaligned\end{equation*}
Then
$$\aligned\langle \Phi'_{\epsilon_n}(w_{n,1}),\varphi\rangle&=\langle \tilde{\Phi}'_{\epsilon_n}(w_{n}),\varphi(x+x_n)\rangle\\
&=\langle \tilde{\Phi}'_{\epsilon_n}(\bar{z}_{n}),\varphi(x+x_n)\rangle-\langle \tilde{\Phi}'_{\epsilon_n}(z_0),\varphi(x+x_n)\rangle+o_n(1)
\|\varphi\|_{1,\epsilon_n}\\
&=0-\langle \Phi'_{0}(z_0),\varphi(x+x_n)\rangle+o_n(1)\|\varphi\|_{1,\epsilon_n}
=o_n(1)\|\varphi\|_{1,\epsilon_n}.
\endaligned$$
Therefore,
\begin{equation}\label{4.12} \Phi'_{\epsilon_n}(w_{n,1})\rightarrow0.\end{equation}
Similarly
$${\Phi}_{\epsilon_n}({z}_{n})-{\Phi}_{\epsilon_n}(w_{n,1})
-\Phi_{0}(z_0)=\tilde{\Phi}_{\epsilon_n}(\bar{z}_{n})-\tilde{\Phi}_{\epsilon_n}(w_{n})
-\Phi_{0}(z_0)\rightarrow0.$$
By (\ref{4.11}) we get
\begin{equation}\label{4.13}{\Phi}_{\epsilon_n}(w_{n,1})\rightarrow0.\end{equation}
Below we show that $\{w_{n,1}\}$ is nonvanishing. Otherwise, if  $\{w_{n,1}:=(w^1_{n,1},w^2_{n,1})\}$ is vanishing, then $w^1_{n,1}\rightarrow0$ and $w^2_{n,1}\rightarrow0$  in $L^r(\mathbb{R}^2)$ for any $r>2$. As (\ref{5.3}) we have
$\int_{\mathbb{R}^2}\bar{g}(\epsilon x,w^2_{n,1})w^1_{n,1}=o_n(1).$
By (\ref{4.12}) and (\ref{4.13}) we have
$$o_n(1)=\langle {\Phi}'_{\epsilon_n}(w_{n,1}),(w^{1}_{n,1},0)\rangle=\|w^{1}_{n,1}\|^2_{\epsilon_n}
-\int_{\mathbb{R}^2}\bar{g}(\epsilon x,w^2_{n,1})w^1_{n,1}=\|w^{1}_{n,1}\|^2_{\epsilon_n}+o_n(1).
$$
Then
$w^{1}_{n,1}\rightarrow0$ in $H^1(\mathbb{R}^2)$. Similarly, $w^{2}_{n,1}\rightarrow0$ in $H^1(\mathbb{R}^2)$,
contradicts with the hypothesis that $w_n\not\rightarrow(0,0)$ in $E$. Hence, there exist $\{x^2_n\}\subset\mathbb{R}^2$ and $\delta>0$ such that
\begin{equation}\label{4.15}\int_{B_1(x^2_n)}[(w^1_{n,1})^2+(w^2_{n,1})^2]
\geq\delta.\end{equation}
Set $w_{n,2}=w_{n,1}(x+x^2_n)$ and assume $w_{n,2}=(w^1_{n,2},w^2_{n,2})\rightharpoonup w_2=(w^1_2,w^2_2)\neq(0,0)$ in $E$. Then $w_{n,2}$ satisfies
\begin{equation}\label{4.14}\tilde{\tilde{\Phi}}'_{\epsilon_n}(w_{n,2})\rightarrow0,\end{equation}
with energy functional
$$\aligned\tilde{\tilde{\Phi}}_{\epsilon_n}(u,v)
=&\int_{\mathbb{R}^2}(\nabla u\nabla v+\tilde{\tilde{V}}_{\epsilon_n}(x)uv)\\&
-\int_{\mathbb{R}^2}\bar{F}(\epsilon_nx+\epsilon_nx^2_n, u)-\int_{\mathbb{R}^2}\bar{G}(\epsilon_nx+\epsilon_nx^2_n, v),\ \ \text{for all}\ (u,v)\in E,\endaligned$$
where $\tilde{\tilde{V}}_{\epsilon_n}(x)=V(\epsilon_n x+\epsilon_nx^2_n)$. According to (\ref{4.15}), we next discuss for two cases.

{\bf Case 1}:  $\int_{B_1(x^2_n)}\bigl[(1-\chi(\epsilon_n x))(w^1_{n,1})^2+(1-\chi(\epsilon_n x))(w^2_{n,1})^2\bigr]\geq{\delta}/{2}$.
In this case, $B_1(x^2_n)\cap \{\mathbb{R}^2\setminus supp\chi(\epsilon_n x) \} \neq\emptyset$.
Then $|\epsilon_n x^2_n|\rightarrow\infty $ or $\epsilon_n x^2_n\rightarrow x'_0\in \mathbb{R}^2\backslash\Lambda$ and further assume $x'_0\not\in \partial\Lambda$ since we can argue as the following Case 2 if $x'_0\in \partial\Lambda$.  Since $V\in L^\infty(\mathbb{R}^2)$, we may assume $V_2=\lim_{n\rightarrow\infty}V(\epsilon_n x^2_n)$. By (\ref{4.14}) we infer
$w_2=(w^1_2,w^2_2)$ satisfies
\begin{equation*}\aligned
\left\{ \begin{array}{lll}
-\Delta w^1_2+V_2w^1_2=\tilde{g}(w^2_2)\ & \text{in}\quad \mathbb{R}^2,\\
-\Delta w^2_2+V_2w^2_2=\tilde{f}(w^1_2)\ & \text{in}\quad \mathbb{R}^2.
\end{array}\right.\endaligned
\end{equation*}
As before we have $w^1_2=w^2_2=0$, contradicts $w_2=(w^1_2,w^2_2)\neq(0,0)$.

{\bf Case 2}:  $\int_{B_1(x^2_n)}\bigl[\chi(\epsilon_n x)(w^1_{n,1})^2+\chi(\epsilon_n x)(w^2_{n,1})^2\bigr]\geq{\delta}/{2}$.
If this case occurs, we may assume $\epsilon_nx^2_n\rightarrow x'_0\in \bar{\Lambda}$. Using (\ref{4.14}) we deduce that
$w_2=(w^1_2,w^2_2)$ satisfies
\begin{equation}\label{4.16}\aligned
\left\{ \begin{array}{lll}
-\Delta w^1_2+V(x'_0)w^1_2={g}_{1,\infty}(x,w^2_2)\ & \text{in}\quad \mathbb{R}^2,\\
-\Delta w^2_2+V(x'_0)w^2_2={f}_{1,\infty}(x,w^1_2)\ & \text{in}\quad \mathbb{R}^2.
\end{array}\right.\endaligned
\end{equation}
where ${g}_{1,\infty}(x,s)=\chi_{1,\infty} g(s)+(1-\chi_{1,\infty})\tilde{g}(s)$, ${f}_{1,\infty}(x,s)=\chi_{1,\infty} f(s)+(1-\chi_{1,\infty})\tilde{f}(s)$ and $\chi_{1,\infty}$ is either a characteristic
function of a half-space of $\mathbb{R}^2$
or $\chi_{1,\infty}\equiv1$. Denote
$\Phi_{1,\infty}$ to be the associate energy functional of (\ref{4.16})
$$\Phi_{1,\infty}(u,v)=\int_{\mathbb{R}^2}(\nabla u\nabla v+V(x'_0)uv)-\int_{\mathbb{R}^2}[F_{1,\infty}(x,u)+
G_{1,\infty}(x,v)],\quad\forall (u,v)\in E,$$
where $F_{1,\infty}(x,s)=\int^s_0f_{1,\infty}(x,t)dt$, $G_{1,\infty}(x,s)=\int^s_0g_{1,\infty}(x,t)dt$.
Observe that $$\int_{\mathbb{R}^2}F_{1,\infty}(x,u)\leq\int_{\mathbb{R}^2}F(u),\quad \int_{\mathbb{R}^2}G_{1,\infty}(x,v)\leq\int_{\mathbb{R}^2}G(v),$$ we have
$$\Phi_{1,\infty}(u,v)=\tilde{\Phi}_{V(x'_0)}(u,v)+\frac{V(x'_0)}{V_0}
\int_{\mathbb{R}^2}(G(v)+F(u))-
\int_{\mathbb{R}^2}[F_{1,\infty}(x,u)
+G_{1,\infty}(x,v)],$$
{for all} $(u,v)\in E$, where $\tilde{\Phi}_{V(x'_0)}$ is given in (\ref{bu3}) with $x_0$ replaced by $x'_0$.
As before, define $h_{1,\infty}: H^1(\mathbb{R}^2)\rightarrow H^1(\mathbb{R}^2)$ and $J_{1,\infty}: H^1(\mathbb{R}^2)\rightarrow\mathbb{R}$ by
$$J_{1,\infty}(u)=\Phi_{1,\infty}(u+h_{1,\infty}(u),u-h_{1,\infty}(u))=\max_{v\in H^1(\mathbb{R}^2)}\Phi_{1,\infty}(u+v,u-v).$$
Furthermore, as before we can show that, if $u\in H^1(\mathbb{R}^2)\backslash\{0\}$ satisfies $\langle J'_{1,\infty}(u),u\rangle=0$, then $J''_{1,\infty}(u)[u,u]<0$. Since  $(w^{1}_2,w^{2}_2)\neq(0,0)$ is a critical point of $\Phi_{1,\infty}$, we then infer that $\frac{w^{1}_2+w^{2}_2}{2}$ is a critical point of $J_{1,\infty}$ and $J_{1,\infty}(\frac{w^{1}_2+w^{2}_2}{2})=\max_{t>0}J_{1,\infty}(t\frac{w^{1}_2+w^{2}_2}{2})$. Similar to (\ref{4.2.5}), no matter $f$ is  superlinear or asymptotically linear, there always exists $t_1>0$ such that $\tilde{J}_{V(x'_0)}(t_1 \frac{w^{1}_2+w^{2}_2}{2})=\max_{t>0}\tilde{J}_{V(x'_0)}(t \frac{w^{1}_2+w^{2}_2}{2})$, where $\tilde{J}_{V(x'_0)}$ is given in (\ref{3.7.1}) with $x_0$ replaced by $x'_0$.
Define $$\psi_1=t_1\frac{w^{1}_2+w^{2}_2}{2}+\tilde{h}_{V(x'_0)}(t_1\frac{w^{1}_2+w^{2}_2}{2}),
\quad \psi_2=t_1\frac{w^{1}_2+w^{2}_2}{2}-\tilde{h}_{V(x'_0)}(t_1
\frac{w^{1}_2+w^{2}_2}{2}),$$
where $\tilde{h}_{V(x'_0)}$ is given in (\ref{3.7.3}) with $x_0$ replaced by $x'_0$.
In view of (\ref{3.7.2}) we deduce
\begin{equation}\label{4.4}\aligned &\Phi_{1,\infty}(w^{1}_2,w^{2}_2)=J_{1,\infty}(\frac{w^{1}_2+w^{2}_2}{2})\geq J_{1,\infty}(t_1\frac{w^{1}_2+w^{2}_2}{2})\\
=&\Phi_{1,\infty}\Bigl(t_1\frac{w^{1}_2+w^{2}_2}{2}+{h}_{1,\infty}(t_1\frac{w^{1}_2+w^{2}_2}{2}),
t_1\frac{w^{1}_2+w^{2}_2}{2}-{h}_{1,\infty}(t_1\frac{w^{1}_2+w^{2}_2}{2})\Bigr)
\\ \geq& \Phi_{1,\infty}(\psi_1,
\psi_2)= \tilde{\Phi}_{V(x'_0)}(\psi_1,
\psi_2)
+\int_{\mathbb{R}^2}\frac{V(x'_0)}{V_0}\bigl(G(\psi_2)+F(\psi_1)\bigr)
\\&-\int_{\mathbb{R}^2}G_{1,\infty}(x,\psi_2)
-\int_{\mathbb{R}^2}F_{1,\infty}(x,\psi_1)\\
\geq& \tilde{\Phi}_{V(x'_0)}(\psi_1,
\psi_2)=\tilde{J}_{V(x'_0)}(t_1\frac{w^{1}_2+w^{2}_2}{2})\geq c_{V_0}.\endaligned\end{equation}
On the other hand, denote $\hat{F}_{1,\infty}( x,s)=\frac12f_{1,\infty}(x,s)s-F_{1,\infty}(x,s)$ and
$\hat{G}_{1,\infty}( x,s)=\frac12g_{1,\infty}(x,s)s-G_{1,\infty}(x,s)$ for all $(x,s)\in\mathbb{R}^2\times \mathbb{R}$, from (\ref{4.12}), (\ref{4.13}) and Fatou Lemma it follows that
\begin{equation*}\aligned &o_n(1)={\Phi}_{\epsilon_n}(w_{n,1})-\frac12\langle {\Phi}'_{\epsilon_n}(w_{n,1}),w_{n,1}\rangle=
\int_{\mathbb{R}^2}(\hat{F}(\epsilon_n x,w^1_{n,1})+\hat{G}(\epsilon_n x,w^2_{n,1}))\\ =&\int_{\mathbb{R}^2}(\hat{F}(\epsilon_n x+\epsilon_n x^2_n,w^1_{n,2})+\hat{G}(\epsilon_n x+\epsilon_n x^2_n,w^2_{n,2}))\\
\geq& \int_{\mathbb{R}^2}(\hat{F}_{1,\infty}( x,w^1_2)+\hat{G}_{1,\infty}( x,w^2_2))=\Phi_{1,\infty}(w_2)-\frac12\langle \Phi'_{1,\infty}(w_2),w_2\rangle=\Phi_{1,\infty}(w_2),\endaligned\end{equation*}
contradicts with (\ref{4.4}) and so $\bar{z}_n\rightarrow z_0$ in $E$.\ \ \ \ \ \ $\Box$

\begin{lemma}\label{l4.3} Let ${z}_{n}=(u_n,v_n)$ are nontrivial solutions of (\ref{2.3}) obtained in Lemma \ref{l3.5.0} with
$\epsilon_n\rightarrow0$, and $\bar{{z}}_n=z_n(\cdot+x_n)\rightarrow z_0=(u_0,v_0)$ in $E$, where $x_n$ and $(u_0,v_0)$ are obtained in Lemma \ref{l4.2}. Then
$$u_n(x+x_n)\rightarrow0,\ v_n(x+x_n)\rightarrow0\ \text{ uniformly in } n\ \text{as}\ |x|\rightarrow\infty.$$
In addition
$$\sup_{n\geq1}(|u_n|_\infty+|v_n|_\infty)<+\infty.$$
\end{lemma}
{\bf Proof}: Let $\bar{u}_n(\cdot)=u_n(\cdot+x_n)$ and $\bar{v}_n(\cdot)=v_n(\cdot+x_n)$. Note that $\bar{u}_n$ is a weak solution of the following problem
\begin{equation*}-\Delta U+V(\epsilon_nx+\epsilon_nx_n)U=\bar{g}(\epsilon_nx+\epsilon_nx_n, \bar{v}_n) \ \text{in}\ B_2,\ U-\bar{u}_n\in H^1_0(B_2),\end{equation*}
where $B_2=B_2(0)$. Moreover, for any $p\geq2$, we have
\begin{equation}\label{4.0.1}\aligned
\|\bar{u}_n\|_{W^{2,p}(B_1)}&\leq C(\|\bar{g}(\epsilon_nx+\epsilon_nx_n, \bar{v}_n)\|_{L^p(B_2)}+\|\bar{u}_n\|_{L^p(B_2)})\\&\leq C(\|g( \bar{v}_n)\|_{L^p(B_2)}+\|\bar{u}_n\|_{L^p(B_2)}).\endaligned\end{equation}
Denote $B_1=B_1(0)$. By the Sobolev embedding theorem, if $p>2$, we get $\bar{u}_n\in C^{1,\gamma}(\bar{B}_1)$ for some $\gamma\in (0,1)$ and there exists $C$ (independent of $n$) such that
 $\|\bar{u}_n\|_{C^{1,\gamma}(\bar{B}_1)}\leq C\|\bar{u}_n\|_{W^{2,p}(B_1)}.$
Using (\ref{4.0.1})  we get
\begin{equation*}
\|\bar{u}_n\|_{C^{1,\gamma}(\bar{B}_1)}\leq C(|g( \bar{v}_n)|_{p}+|\bar{u}_n|_{p}).\end{equation*}
Assume $\|s\bar{v}_n+(1-s)v_0\|^2_{\epsilon_n}\leq M$, where $s\in (0,1)$.
For some $\alpha\in (0,\frac{2\pi}{pM})$, the Trudinger-Moser inequality implies that
$$\int_{\mathbb{R}^2}\Bigl[e^{2\alpha p(s\bar{v}_n+(1-s)v_0)^2}-1\Bigr]\leq
\int_{\mathbb{R}^2}\bigl[e^{2\alpha pM\frac{(s\bar{v}_n+(1-s)v_0)^2}{\|s\bar{v}_n+(1-s)v_0\|^2_{\epsilon_n}}}-1\bigr]\leq C.$$
As (\ref{2.2.0}), for the above $\alpha$, for any $\delta>0$ there exists $C_\delta>0$ such that
$$\aligned
\int_{\mathbb{R}^2}|g(\bar{v}_n)-g(v_0)|^p
&=\int_{\mathbb{R}^2}|g'(s\bar{v}_n+(1-s)v_0)|^p|\bar{v}_n-v_0|^p\\
&\leq\delta|\bar{v}_n-v_0|^p_p+C_\delta\Bigl(\int_{\mathbb{R}^2}[e^{2\alpha p(s\bar{v}_n+(1-s)v_0)^2}-1]\Bigr)^{\frac12}\|\bar{v}_n-v_0|^{p}_{2p} \rightarrow0.\endaligned$$
Recalling that $\bar{v}_n\rightarrow v_0$ in $H^1(\mathbb{R}^2)$, we get
\begin{equation}\label{4.2.2}|g(\bar{v}_n)|^p_p\rightarrow|g(v_0)|^p_p.\end{equation} Hence
\begin{equation}\label{4.1.4}
\sup_{n\geq1}\|\bar{u}_n\|_{C^{1,\gamma}(\bar{B}_1)}<+\infty.
\end{equation}
Next we prove that
\begin{equation}\label{4.1.5}\bar{u}_n(x)\rightarrow0 \text{\ uniformly as}\ |x|\rightarrow\infty.\end{equation}
  It suffices to prove that for any
$\delta>0$, there exists $R>0$ such that $|\bar{u}_n(x)|\leq \delta$,\ $\forall n\geq1$, $|x|>R$. Otherwise, there exists $\{P_n\}\subset\mathbb{R}^2$ such that $|P_n|\rightarrow\infty$ and
$\liminf_{n\rightarrow\infty}|\bar{u}_n(P_n)|>0$. Let
$\tilde{u}_n(x)=\bar{u}_n(x+P_n)$ and $\tilde{v}_n(x)=\bar{v}_n(x+P_n)$, then
\begin{equation*} -\Delta \tilde{u}_n+V(\epsilon_nx+\epsilon_nx_n+\epsilon_nP_n)\tilde{u}_n=
\bar{g}(\epsilon_nx+\epsilon_nx_n+\epsilon_nP_n,\tilde{v}_n),\ \quad \tilde{u}_n\in H^1(\mathbb{R}^2).\end{equation*}
 Suppose that $\tilde{u}_n\rightharpoonup \tilde{u}_0$ in $H^1(\mathbb{R}^2)$ and claim $\tilde{u}_0\not\equiv0$. Indeed, for any $n\geq1$, $\tilde{{u}}_n$ is a weak solution of the problem
\begin{equation*}-\Delta U+V(\epsilon_nx+\epsilon_nx_n+\epsilon_nP_n)U=\bar{g}(\epsilon_nx+
\epsilon_nx_n+\epsilon_nP_n, \tilde{{v}}_n) \ \text{in}\ B_2,\ U-\tilde{u}_n\in H^1_0(B_2).\end{equation*}
In the same way as (\ref{4.1.4}) we have
$\sup_{n\geq1}\|\tilde{{u}}_n\|_{C^{1,\gamma}(\bar{B}_1)}<+\infty.$
Up to a subsequence, $\tilde{{u}}_n\rightarrow \tilde{{u}}_0$ uniformly in $\bar{B}_1$. Thus,
$$\tilde{u}_{0}(0)=\liminf_{n\rightarrow\infty}\tilde{u}_{n}(0)
=\liminf_{n\rightarrow\infty}\bar{{u}}_{n}(P_n)\neq0.$$
On the other hand, note that $\bar{u}_n\rightarrow u_0$ in $H^1(\mathbb{R}^2)$, for any fixed $R>0$ and $n$ large enough, we infer
$$\aligned
o_n(1)+\int_{\mathbb{R}^2}u^2_0&=\int_{\mathbb{R}^2}\bar{u}^2_n
\geq\int_{B_R(0)}(\bar{u}^2_n+\tilde{u}^2_n)
=\int_{B_R(0)}({u}^2_0+\tilde{{u}}^2_0)+o_n(1).\endaligned$$
Since $R$ is arbitrary, we know $\tilde{u}_0\equiv0$, contradicts $\tilde{u}_{0}(0)\neq0$ before. Thus, (\ref{4.1.5}) holds true. Combining with (\ref{4.1.4}) we infer that
$\sup_{n\geq1}|u_n|_\infty=\sup_{n\geq1}|\bar{u}_n|_\infty<+\infty$. Similarly, $\sup_{n\geq1}|v_n|_\infty<+\infty$.\ \ \ \ \ $\Box$

\begin{lemma}\label{l4.4} There exist $\delta,M>0$ independent of $n$ such that $\delta\leq\min\{|u_n|_\infty,|v_n|_\infty\}<|u_n|_\infty+|v_n|_\infty\leq M$.\end{lemma}
{\bf Proof}: In view of Lemma \ref{l4.3}, it suffices to show that
$\min\{|u_n|_\infty,|v_n|_\infty\}\geq\delta$. We argue by contradiction and assume
$\min\{|u_n|_\infty,|v_n|_\infty\}\rightarrow0.$
Without loss of generality, we assume that $|v_n|_\infty\rightarrow0$. By (H$_1$) we get
$$\int_{\mathbb{R}^2}(|\nabla u_n|^2+V(\epsilon_n x)u^2_n)=\int_{\mathbb{R}^2}\bar{g}(\epsilon_n x,v_n)u_n\leq o_n(1)|u_n|_2|v_n|_2.$$
Hence $\|u_n\|_\epsilon\rightarrow0$. Moreover, by (\ref{2.1.0}) with $q=1$ and $\alpha<\frac{4\pi}{3}$ we have
\begin{equation*}\aligned \int_{\mathbb{R}^2}|\nabla v_n|^2+V(\epsilon_n x)v^2_n&=\int_{\mathbb{R}^2}\bar{f}(\epsilon_n x,u_n)v_n\\& \leq \delta|u_n|_2|v_n|_2+C_\delta|u_n|_3|v_n|_3
\int_{\mathbb{R}^2}[e^{3\alpha u^2_n}-1]\rightarrow0.\endaligned\end{equation*}
Then
$v_n\rightarrow0$ in $H^1(\mathbb{R}^2)$ and so
$c_{\epsilon_n}=\Phi_{\epsilon_n}(u_n,v_n)\rightarrow0,$
contradicts $c_{\epsilon_n}\geq\tau>0$.\ \ \ \ \ \ \ $\Box$

By Lemma \ref{l4.4}, we may assume $\{y_n\}\subset\mathbb{R}^2$ satisfies
\begin{equation}\label{5.4}|u_n(y_n)|+|v_n(y_n)|=
\max_{\mathbb{R}^2}(|u_n(x)|+|v_n(x)|).\end{equation} Then the following result holds.
\begin{lemma}\label{l4.5} (1) $\lim_{n\rightarrow\infty} dist(\epsilon_n y_n,\mathcal{V})=0$, and $z_n(\cdot+y_n)\rightarrow \hat{z}_0$ in $E$ as $n\rightarrow\infty$, where $\hat{z}_0$ is a ground state of the limit equation (\ref{1.1.1});\\
\noindent (2) $u_n(x+y_n)\rightarrow0$ and $v_n(x+y_n)\rightarrow0$ uniformly in $n$ as $|x|\rightarrow\infty$.
\end{lemma}
{\bf Proof}: (1)  Firstly we claim that there exist $\mu, R_1>0$ such that
\begin{equation}\label{4.2.0} \lim_{n\rightarrow\infty}\int_{B_{R_1}(y_n)}(u^2_n+v^2_n)\geq\mu.\end{equation}
Argue by contradiction we assume that, for any $R>0$,
$$\lim_{n\rightarrow\infty}\int_{B_{R}(y_n)}(u^2_n+v^2_n)=0.$$
Let $\hat{u}_n(\cdot)=u_n(\cdot+y_n)$ and $\hat{v}_n(\cdot)=v_n(\cdot+y_n)$, then $\hat{u}_n, \hat{v}_n\rightarrow0$ in $L^2_{loc}(\mathbb{R}^2)$ as $n\rightarrow\infty$. Observe that $\hat{u}_n$ is a weak solution of the following problem
\begin{equation*}-\Delta U+V(\epsilon_nx+\epsilon_ny_n)U=\bar{g}(\epsilon_nx+\epsilon_ny_n, \hat{{v}}_n) \ \text{in}\ B_2,\ U-\hat{u}_n\in H^1_0(B_2).\end{equation*}
As in Lemma \ref{l4.3}, by standard elliptic regularity we get $\hat{u}_n\in C^{1,\gamma}(\bar{B}_1)$ for some $\gamma\in (0,1)$ and there exists $C$ (independent of $n$) such that
\begin{equation*}\|\hat{u}_n\|_{C^{1,\gamma}(\bar{B}_1)}\leq C(|g( \hat{v}_n)|_{p}+|\hat{u}_n|_{p}).\end{equation*}
Using Lemma \ref{l4.2}, $\bar{z}_n=z_n(x+x_n)\rightarrow z_0$ in $E$. Moreover, by (\ref{4.2.2}) we get
$|g(\hat{v}_n)|_{p}\rightarrow|g(v_0)|_p$. Then
$ \sup_{n\geq1}\|\hat{u}_n\|_{C^{1,\gamma}(\bar{B}_1)}<+\infty.$
Since $\hat{u}_n\rightarrow0$ in $L^2(B_1)$, we get ${\hat{u}_n}\rightarrow0$ uniformly in $B_1$. In particular, $\hat{u}_n(0)=u_n(y_n)\rightarrow0$. Similarly, $\hat{v}_n(0)=v_n(y_n)\rightarrow0$. Then we obtain
$$\lim_{n\rightarrow\infty}\max_{x\in\mathbb{R}^2}(|u_n(x)|+|v_n(x)|)
=\lim_{n\rightarrow\infty}(|u_n(y_n)|+|v_n(y_n)|)=0,$$
contradicts with Lemma \ref{l4.4}. Thus (\ref{4.2.0}) holds true.
Then we may assume $z_n(\cdot+y_n)=(\hat{u}_n(\cdot),\hat{v}_n(\cdot))\rightharpoonup \hat{z}_0\neq(0,0)$. Arguing as Claim 1 and 2 in the proof of Lemma \ref{l4.2},  we deduce that $\epsilon_ny_n\rightarrow y_0\in \Lambda$, $V(y_0)=V_0$ and $z_n(\cdot+y_n)\rightarrow \hat{z}_0$ in $E$, where $\hat{z}_0$ is a ground state of the limit equation (\ref{1.1.1}).

(2) Taking similar arguments as in the proof of (\ref{4.1.5}), we can infer desired results.\ \ \ \ \ $\Box$

\begin{lemma} \label{l5.4} Let $z_n=(u_n,v_n)$ be nontrivial solutions of (\ref{2.3}) obtained in Lemma \ref{l3.5.0} with $\epsilon_n\rightarrow0$, then $u_n$ and $v_n$ have maximum points  $x^1_n$ and $x^2_n$ respectively, and $$\lim_{n\rightarrow\infty}dist(\epsilon_nx^i_n,\mathcal{V})=0,\quad i=1,2.$$
\end{lemma}
{\bf Proof}: Note that $u_n(x+x_n)\rightarrow0$ and $v_n(x+x_n)\rightarrow0$ as $|x|\rightarrow\infty$ and the fact that
$\min\{|u_n|_\infty,|v_n|_\infty\}\geq\delta>0$, then $u_n(\cdot+x_n)$ and $v_n(\cdot+x_n)$ has maximum points $P^1_n$ and $P^2_n$ respectively and $|P^1_n|,|P^2_n|\leq R_0$ for some $R_0>0$. Then $u_n$ and $v_n$ has maximum points $x^1_n$ and $x^2_n$ respectively. In addition, $x^i_n=P^i_n+x_n$, $i=1,2$. Since $dist(\epsilon_nx_n, \mathcal{V})=0$, we get $dist(\epsilon_nx^i_n, \mathcal{V})=o_n(1)$. \ \ \ \ \ $\Box$

{\bf Proof of Theorem 1.1} Firstly we need to show that the solutions $z_n=(u_n,v_n)$ of the penalized problem (\ref{2.3}) obtained in Lemma \ref{l3.5.0} are actually the solutions of the problem (\ref{2.1}). We claim that there exists $\widetilde{\epsilon}>0$ such that for any $\epsilon\in (0,\widetilde{\epsilon})$ and any solution $(u,v)$ of system (\ref{2.3}) there holds
\begin{equation}\label{4.3.2}\|u\|_{L^\infty(\mathbb{R}^2
\backslash{\Lambda_\epsilon})}<a_1,\quad \|v\|_{L^\infty(\mathbb{R}^2
\backslash{\Lambda_\epsilon})}<a_2.\end{equation}
Argue by contradiction assume there exist $\epsilon_n\rightarrow0$, $(u_{n},v_n)$ satisfies $\Phi'_{\epsilon_n}(u_n,v_n)=0$ and without loss of generality we further suppose
\begin{equation}\label{4.3.3}\|u_n\|_{{L^\infty(\mathbb{R}^2
\backslash{\Lambda_{\epsilon_n}})}}\geq a_1.\end{equation}
As in Lemma \ref{l4.2}, there exists $\{x'_n\}\subset\mathbb{R}^2$  such that $\epsilon_nx'_n\rightarrow y_0\in \mathcal{V}$ and $(u_n(\cdot+x'_n),v_n(\cdot+x'_n)\rightarrow(u_0,v_0)$ in $E$.

Since the set $\Lambda$ is open, choose $r>0$ such that $B_r(y_0)\subset B_{2r}(y_0)\subset\Lambda$, we have
$$B_{\frac{2r}{\epsilon_n}}(\frac{y_0}{\epsilon_n})
=\frac{1}{\epsilon_n}B_{2r}(y_0)\subset \Lambda_{\epsilon_n}.$$
Moreover, for any $z\in B_{\frac{r}{\epsilon_n}}(x'_n)$, we get
$$\Bigl|z-\frac{y_0}{\epsilon_n}\Bigr|\leq|z-x'_n|+
\Bigl|x'_n-\frac{y_0}{\epsilon_n}\Bigr|
<\frac{1}{\epsilon_n}(r+o_n(1))<\frac{2r}{\epsilon_n},$$
for large $n$. Then $B_{\frac{r}{\epsilon_n}}(x'_n)\subset\Lambda_{\epsilon_n}$ for large $n$.
 Arguing as Lemma \ref{l4.3}, $|u_n(x+x'_n)|\rightarrow0$ and $|v_n(x+x'_n)|\rightarrow0$ uniformly in $x$ as $n\rightarrow\infty$. Then
there exists $R>0$ such that $\|u_n\|_{L^\infty({\mathbb{R}^2\backslash{B_{R}(x'_n)}})}<a_1$.
So
for any $n\geq n_0$ satisfying $\frac{r}{\epsilon_n}>R$, we deduce
$$\|u_n\|_{L^\infty({\mathbb{R}^2\backslash{\Lambda_{\epsilon_n}}})}
\leq\|u_n\|_{L^\infty({\mathbb{R}^2\backslash{B_{\frac r{\epsilon_n}}(x'_n)}})}\leq
\|u_n\|_{L^\infty({\mathbb{R}^2\backslash{B_{R}(x'_n)}})}<a_1,$$
which contradicts with (\ref{4.3.3}) and then (\ref{4.3.2}) yields. From the definition of $\bar{f}$ we know $\bar{f}(\epsilon x,u)=f(u)$. Similarly, $\bar{g}(\epsilon x,v)=g(v)$. Then $(u,v)$ is also a solution of the problem (\ref{2.1}). Hence $(\varphi_\epsilon(x),\psi_\epsilon(x))=(u(\frac{x}{\epsilon}),
v(\frac{x}{\epsilon}))$ is a solution of the system (\ref{1.1}), and Lemmas \ref{l4.2}-\ref{l5.4} imply that the conclusions of Theorem 1.1 hold true.\ \ \ \ \ $\Box$

{\bf Proof of Theorem 1.2} Assume that $z_n=(u_n,v_n)$ are nontrivial solutions of (\ref{2.1}) obtained in Theorem 1.1 with $\epsilon_n\rightarrow0$, $x^1_n$ and $x^2_n$ are maximum points of $u_n$ and $v_n$ respectively given in Lemma \ref{l5.4}. Using Lemma \ref{l5.3} and taking similar arguments as in the proof of Propositions 3.13 and 3.14 in \cite{ZhangCPDE}, we can deduce that for $\epsilon_n>0$ small enough, $u_nv_n>0$ in $\mathbb{R}^2$,
 $\lim_{n\rightarrow\infty}|x^1_n-x^2_n|=0,$
and for some $c,C>0$,
 $$|{u}_{n}(x)|\leq C e^{-c|x-x^1_n|},\quad \ |{v}_{n}(x)|\leq C e^{-c|x-x^2_n|},\ \ \forall x\in\mathbb{R}^2.$$ Moreover, $x^1_n$ and $x^2_n$ are unique.
Letting $(\varphi_\epsilon(x),\psi_\epsilon(x))=(u(\frac{x}{\epsilon}),
v(\frac{x}{\epsilon}))$,  the conclusions of Theorem 1.2 yield. \ \ \ \ \ \ $\Box$

\end{document}